\input amstex
%

\def\next{AMS-SEKR}\ifx\styname\next \endinput\fi
\catcode`\@=11
\def\styname{AMS-SEKR}
\def\styversion{2.0}
{\W@{}\W@{\styname.STY - Version \styversion}\W@{}}
\hyphenation{acad-e-my acad-e-mies af-ter-thought anom-aly anom-alies
an-ti-deriv-a-tive an-tin-o-my an-tin-o-mies apoth-e-o-ses apoth-e-o-sis
ap-pen-dix ar-che-typ-al as-sign-a-ble as-sist-ant-ship as-ymp-tot-ic
asyn-chro-nous at-trib-uted at-trib-ut-able bank-rupt bank-rupt-cy
bi-dif-fer-en-tial blue-print busier busiest cat-a-stroph-ic
cat-a-stroph-i-cally con-gress cross-hatched data-base de-fin-i-tive
de-riv-a-tive dis-trib-ute dri-ver dri-vers eco-nom-ics econ-o-mist
elit-ist equi-vari-ant ex-quis-ite ex-tra-or-di-nary flow-chart
for-mi-da-ble forth-right friv-o-lous ge-o-des-ic ge-o-det-ic geo-met-ric
griev-ance griev-ous griev-ous-ly hexa-dec-i-mal ho-lo-no-my ho-mo-thetic
ideals idio-syn-crasy in-fin-ite-ly in-fin-i-tes-i-mal ir-rev-o-ca-ble
key-stroke lam-en-ta-ble light-weight mal-a-prop-ism man-u-script
mar-gin-al meta-bol-ic me-tab-o-lism meta-lan-guage me-trop-o-lis
met-ro-pol-i-tan mi-nut-est mol-e-cule mono-chrome mono-pole mo-nop-oly
mono-spline mo-not-o-nous mul-ti-fac-eted mul-ti-plic-able non-euclid-ean
non-iso-mor-phic non-smooth par-a-digm par-a-bol-ic pa-rab-o-loid
pa-ram-e-trize para-mount pen-ta-gon phe-nom-e-non post-script pre-am-ble
pro-ce-dur-al pro-hib-i-tive pro-hib-i-tive-ly pseu-do-dif-fer-en-tial
pseu-do-fi-nite pseu-do-nym qua-drat-ics quad-ra-ture qua-si-smooth
qua-si-sta-tion-ary qua-si-tri-an-gu-lar quin-tes-sence quin-tes-sen-tial
re-arrange-ment rec-tan-gle ret-ri-bu-tion retro-fit retro-fit-ted
right-eous right-eous-ness ro-bot ro-bot-ics sched-ul-ing se-mes-ter
semi-def-i-nite semi-ho-mo-thet-ic set-up se-vere-ly side-step sov-er-eign
spe-cious spher-oid spher-oid-al star-tling star-tling-ly
sta-tis-tics sto-chas-tic straight-est strange-ness strat-a-gem strong-hold
sum-ma-ble symp-to-matic syn-chro-nous topo-graph-i-cal tra-vers-a-ble
tra-ver-sal tra-ver-sals treach-ery turn-around un-at-tached un-err-ing-ly
white-space wide-spread wing-spread wretch-ed wretch-ed-ly Brown-ian
Eng-lish Euler-ian Feb-ru-ary Gauss-ian Grothen-dieck Hamil-ton-ian
Her-mit-ian Jan-u-ary Japan-ese Kor-te-weg Le-gendre Lip-schitz
Lip-schitz-ian Mar-kov-ian Noe-ther-ian No-vem-ber Rie-mann-ian
Schwarz-schild Sep-tem-ber
form per-iods Uni-ver-si-ty cri-ti-sism for-ma-lism}
\Invalid@\nofrills
\Invalid@\usualspace
\newif\ifnofrills@
\def\nofrills@#1#2{\relaxnext@
  \DN@{\ifx\next\nofrills
    \nofrills@true\let#2\relax\DN@\nofrills{\nextii@}%
  \else
    \nofrills@false\def#2{#1}\let\next@\nextii@\fi
\next@}}
\def\usualspace@#1{\ifnofrills@\def\usualspace{#1}\fi}
\def\addto#1#2{\csname \expandafter\eat@\string#1@\endcsname
  \expandafter{\the\csname \expandafter\eat@\string#1@\endcsname#2}}
\newdimen\bigsize@
\def\big@#1#2{{\hbox{$\left#2\vcenter to#1\bigsize@{}%
  \right.\nulldelimiterspace\z@\m@th$}}}
\def\big{\big@\@ne}
\def\Big{\big@{1.5}}
\def\bigg{\big@\tw@}
\def\Bigg{\big@{2.5}}
\def\raggedcenter@{\leftskip\z@ plus.4\hsize \rightskip\leftskip
 \parfillskip\z@ \parindent\z@ \spaceskip.3333em \xspaceskip.5em
 \pretolerance9999\tolerance9999 \exhyphenpenalty\@M
 \hyphenpenalty\@M \let\\\linebreak}
\def\upperspecialchars{\def\ss{SS}\let\i=I\let\j=J\let\ae\AE\let\oe\OE
  \let\o\O\let\aa\AA\let\l\L}
\def\uppercasetext@#1{%
  {\spaceskip1.2\fontdimen2\the\font plus1.2\fontdimen3\the\font
   \upperspecialchars\uctext@#1$\m@th\aftergroup\eat@$}}
\def\uctext@#1$#2${\endash@#1-\endash@$#2$\uctext@}
\def\endash@#1-#2\endash@{\uppercase{#1}\if\notempty{#2}--\endash@#2\endash@\fi}
\def\runaway@#1{\DN@{#1}\ifx\envir@\next@
  \Err@{You seem to have a missing or misspelled \string\end#1 ...}%
  \let\envir@\empty\fi}
\newif\iftemp@
\def\notempty#1{TT\fi\def\test@{#1}\ifx\test@\empty\temp@false
  \else\temp@true\fi \iftemp@}
\font@\tensmc=cmcsc10
\font@\sevenex=cmex7
\font@\sevenit=cmti7
\font@\eightrm=cmr8 
\font@\sixrm=cmr6 
\font@\eighti=cmmi8     \skewchar\eighti='177 
\font@\sixi=cmmi6       \skewchar\sixi='177   
\font@\eightsy=cmsy8    \skewchar\eightsy='60 
\font@\sixsy=cmsy6      \skewchar\sixsy='60   
\font@\eightex=cmex8
\font@\eightbf=cmbx8 
\font@\sixbf=cmbx6   
\font@\eightit=cmti8 
\font@\eightsl=cmsl8 
\font@\eightsmc=cmcsc8
\font@\eighttt=cmtt8 

\loadmsam
\loadmsbm
\loadeufm
\UseAMSsymbols
\newtoks\tenpoint@
\def\tenpoint{\normalbaselineskip12\p@
 \abovedisplayskip12\p@ plus3\p@ minus9\p@
 \belowdisplayskip\abovedisplayskip
 \abovedisplayshortskip\z@ plus3\p@
 \belowdisplayshortskip7\p@ plus3\p@ minus4\p@
 \textonlyfont@\rm\tenrm \textonlyfont@\it\tenit
 \textonlyfont@\sl\tensl \textonlyfont@\bf\tenbf
 \textonlyfont@\smc\tensmc \textonlyfont@\tt\tentt
 \textonlyfont@\bsmc\tenbsmc
 \ifsyntax@ \def\big##1{{\hbox{$\left##1\right.$}}}%
  \let\Big\big \let\bigg\big \let\Bigg\big
 \else
  \textfont\z@=\tenrm  \scriptfont\z@=\sevenrm  \scriptscriptfont\z@=\fiverm
  \textfont\@ne=\teni  \scriptfont\@ne=\seveni  \scriptscriptfont\@ne=\fivei
  \textfont\tw@=\tensy \scriptfont\tw@=\sevensy \scriptscriptfont\tw@=\fivesy
  \textfont\thr@@=\tenex \scriptfont\thr@@=\sevenex
        \scriptscriptfont\thr@@=\sevenex
  \textfont\itfam=\tenit \scriptfont\itfam=\sevenit
        \scriptscriptfont\itfam=\sevenit
  \textfont\bffam=\tenbf \scriptfont\bffam=\sevenbf
        \scriptscriptfont\bffam=\fivebf
  \setbox\strutbox\hbox{\vrule height8.5\p@ depth3.5\p@ width\z@}%
  \setbox\strutbox@\hbox{\lower.5\normallineskiplimit\vbox{%
        \kern-\normallineskiplimit\copy\strutbox}}%
 \setbox\z@\vbox{\hbox{$($}\kern\z@}\bigsize@=1.2\ht\z@
 \fi
 \normalbaselines\rm\ex@.2326ex\jot3\ex@\the\tenpoint@}
\newtoks\eightpoint@
\def\eightpoint{\normalbaselineskip10\p@
 \abovedisplayskip10\p@ plus2.4\p@ minus7.2\p@
 \belowdisplayskip\abovedisplayskip
 \abovedisplayshortskip\z@ plus2.4\p@
 \belowdisplayshortskip5.6\p@ plus2.4\p@ minus3.2\p@
 \textonlyfont@\rm\eightrm \textonlyfont@\it\eightit
 \textonlyfont@\sl\eightsl \textonlyfont@\bf\eightbf
 \textonlyfont@\smc\eightsmc \textonlyfont@\tt\eighttt
 \textonlyfont@\bsmc\eightbsmc
 \ifsyntax@\def\big##1{{\hbox{$\left##1\right.$}}}%
  \let\Big\big \let\bigg\big \let\Bigg\big
 \else
  \textfont\z@=\eightrm \scriptfont\z@=\sixrm \scriptscriptfont\z@=\fiverm
  \textfont\@ne=\eighti \scriptfont\@ne=\sixi \scriptscriptfont\@ne=\fivei
  \textfont\tw@=\eightsy \scriptfont\tw@=\sixsy \scriptscriptfont\tw@=\fivesy
  \textfont\thr@@=\eightex \scriptfont\thr@@=\sevenex
   \scriptscriptfont\thr@@=\sevenex
  \textfont\itfam=\eightit \scriptfont\itfam=\sevenit
   \scriptscriptfont\itfam=\sevenit
  \textfont\bffam=\eightbf \scriptfont\bffam=\sixbf
   \scriptscriptfont\bffam=\fivebf
 \setbox\strutbox\hbox{\vrule height7\p@ depth3\p@ width\z@}%
 \setbox\strutbox@\hbox{\raise.5\normallineskiplimit\vbox{%
   \kern-\normallineskiplimit\copy\strutbox}}%
 \setbox\z@\vbox{\hbox{$($}\kern\z@}\bigsize@=1.2\ht\z@
 \fi
 \normalbaselines\eightrm\ex@.2326ex\jot3\ex@\the\eightpoint@}

\font@\twelverm=cmr10 scaled\magstep1
\font@\twelveit=cmti10 scaled\magstep1
\font@\twelvesl=cmsl10 scaled\magstep1
\font@\twelvesmc=cmcsc10 scaled\magstep1
\font@\twelvett=cmtt10 scaled\magstep1
\font@\twelvebf=cmbx10 scaled\magstep1
\font@\twelvei=cmmi10 scaled\magstep1
\font@\twelvesy=cmsy10 scaled\magstep1
\font@\twelveex=cmex10 scaled\magstep1
\font@\twelvemsa=msam10 scaled\magstep1
\font@\twelveeufm=eufm10 scaled\magstep1
\font@\twelvemsb=msbm10 scaled\magstep1
\newtoks\twelvepoint@
\def\twelvepoint{\normalbaselineskip15\p@
 \abovedisplayskip15\p@ plus3.6\p@ minus10.8\p@
 \belowdisplayskip\abovedisplayskip
 \abovedisplayshortskip\z@ plus3.6\p@
 \belowdisplayshortskip8.4\p@ plus3.6\p@ minus4.8\p@
 \textonlyfont@\rm\twelverm \textonlyfont@\it\twelveit
 \textonlyfont@\sl\twelvesl \textonlyfont@\bf\twelvebf
 \textonlyfont@\smc\twelvesmc \textonlyfont@\tt\twelvett
 \textonlyfont@\bsmc\twelvebsmc
 \ifsyntax@ \def\big##1{{\hbox{$\left##1\right.$}}}%
  \let\Big\big \let\bigg\big \let\Bigg\big
 \else
  \textfont\z@=\twelverm  \scriptfont\z@=\tenrm  \scriptscriptfont\z@=\sevenrm
  \textfont\@ne=\twelvei  \scriptfont\@ne=\teni  \scriptscriptfont\@ne=\seveni
  \textfont\tw@=\twelvesy \scriptfont\tw@=\tensy \scriptscriptfont\tw@=\sevensy
  \textfont\thr@@=\twelveex \scriptfont\thr@@=\tenex
        \scriptscriptfont\thr@@=\tenex
  \textfont\itfam=\twelveit \scriptfont\itfam=\tenit
        \scriptscriptfont\itfam=\tenit
  \textfont\bffam=\twelvebf \scriptfont\bffam=\tenbf
        \scriptscriptfont\bffam=\sevenbf
  \setbox\strutbox\hbox{\vrule height10.2\p@ depth4.2\p@ width\z@}%
  \setbox\strutbox@\hbox{\lower.6\normallineskiplimit\vbox{%
        \kern-\normallineskiplimit\copy\strutbox}}%
 \setbox\z@\vbox{\hbox{$($}\kern\z@}\bigsize@=1.4\ht\z@
 \fi
 \normalbaselines\rm\ex@.2326ex\jot3.6\ex@\the\twelvepoint@}

\def\headfonts{\twelvepoint\bf}

\font@\fourteenrm=cmr10 scaled\magstep2
\font@\fourteenit=cmti10 scaled\magstep2
\font@\fourteensl=cmsl10 scaled\magstep2
\font@\fourteensmc=cmcsc10 scaled\magstep2
\font@\fourteentt=cmtt10 scaled\magstep2
\font@\fourteenbf=cmbx10 scaled\magstep2
\font@\fourteeni=cmmi10 scaled\magstep2
\font@\fourteensy=cmsy10 scaled\magstep2
\font@\fourteenex=cmex10 scaled\magstep2
\font@\fourteenmsa=msam10 scaled\magstep2
\font@\fourteeneufm=eufm10 scaled\magstep2
\font@\fourteenmsb=msbm10 scaled\magstep2
\newtoks\fourteenpoint@
\def\fourteenpoint{\normalbaselineskip15\p@
 \abovedisplayskip18\p@ plus4.3\p@ minus12.9\p@
 \belowdisplayskip\abovedisplayskip
 \abovedisplayshortskip\z@ plus4.3\p@
 \belowdisplayshortskip10.1\p@ plus4.3\p@ minus5.8\p@
 \textonlyfont@\rm\fourteenrm \textonlyfont@\it\fourteenit
 \textonlyfont@\sl\fourteensl \textonlyfont@\bf\fourteenbf
 \textonlyfont@\smc\fourteensmc \textonlyfont@\tt\fourteentt
 \textonlyfont@\bsmc\fourteenbsmc
 \ifsyntax@ \def\big##1{{\hbox{$\left##1\right.$}}}%
  \let\Big\big \let\bigg\big \let\Bigg\big
 \else
  \textfont\z@=\fourteenrm  \scriptfont\z@=\twelverm  \scriptscriptfont\z@=\tenrm
  \textfont\@ne=\fourteeni  \scriptfont\@ne=\twelvei  \scriptscriptfont\@ne=\teni
  \textfont\tw@=\fourteensy \scriptfont\tw@=\twelvesy \scriptscriptfont\tw@=\tensy
  \textfont\thr@@=\fourteenex \scriptfont\thr@@=\twelveex
        \scriptscriptfont\thr@@=\twelveex
  \textfont\itfam=\fourteenit \scriptfont\itfam=\twelveit
        \scriptscriptfont\itfam=\twelveit
  \textfont\bffam=\fourteenbf \scriptfont\bffam=\twelvebf
        \scriptscriptfont\bffam=\tenbf
  \setbox\strutbox\hbox{\vrule height12.2\p@ depth5\p@ width\z@}%
  \setbox\strutbox@\hbox{\lower.72\normallineskiplimit\vbox{%
        \kern-\normallineskiplimit\copy\strutbox}}%
 \setbox\z@\vbox{\hbox{$($}\kern\z@}\bigsize@=1.7\ht\z@
 \fi
 \normalbaselines\rm\ex@.2326ex\jot4.3\ex@\the\fourteenpoint@}

\def\chapheadfonts{\fourteenpoint\bf}

\font@\seventeenrm=cmr10 scaled\magstep3
\font@\seventeenit=cmti10 scaled\magstep3
\font@\seventeensl=cmsl10 scaled\magstep3
\font@\seventeensmc=cmcsc10 scaled\magstep3
\font@\seventeentt=cmtt10 scaled\magstep3
\font@\seventeenbf=cmbx10 scaled\magstep3
\font@\seventeeni=cmmi10 scaled\magstep3
\font@\seventeensy=cmsy10 scaled\magstep3
\font@\seventeenex=cmex10 scaled\magstep3
\font@\seventeenmsa=msam10 scaled\magstep3
\font@\seventeeneufm=eufm10 scaled\magstep3
\font@\seventeenmsb=msbm10 scaled\magstep3
\newtoks\seventeenpoint@
\def\seventeenpoint{\normalbaselineskip18\p@
 \abovedisplayskip21.6\p@ plus5.2\p@ minus15.4\p@
 \belowdisplayskip\abovedisplayskip
 \abovedisplayshortskip\z@ plus5.2\p@
 \belowdisplayshortskip12.1\p@ plus5.2\p@ minus7\p@
 \textonlyfont@\rm\seventeenrm \textonlyfont@\it\seventeenit
 \textonlyfont@\sl\seventeensl \textonlyfont@\bf\seventeenbf
 \textonlyfont@\smc\seventeensmc \textonlyfont@\tt\seventeentt
 \textonlyfont@\bsmc\seventeenbsmc
 \ifsyntax@ \def\big##1{{\hbox{$\left##1\right.$}}}%
  \let\Big\big \let\bigg\big \let\Bigg\big
 \else
  \textfont\z@=\seventeenrm  \scriptfont\z@=\fourteenrm  \scriptscriptfont\z@=\twelverm
  \textfont\@ne=\seventeeni  \scriptfont\@ne=\fourteeni  \scriptscriptfont\@ne=\twelvei
  \textfont\tw@=\seventeensy \scriptfont\tw@=\fourteensy \scriptscriptfont\tw@=\twelvesy
  \textfont\thr@@=\seventeenex \scriptfont\thr@@=\fourteenex
        \scriptscriptfont\thr@@=\fourteenex
  \textfont\itfam=\seventeenit \scriptfont\itfam=\fourteenit
        \scriptscriptfont\itfam=\fourteenit
  \textfont\bffam=\seventeenbf \scriptfont\bffam=\fourteenbf
        \scriptscriptfont\bffam=\twelvebf
  \setbox\strutbox\hbox{\vrule height14.6\p@ depth6\p@ width\z@}%
  \setbox\strutbox@\hbox{\lower.86\normallineskiplimit\vbox{%
        \kern-\normallineskiplimit\copy\strutbox}}%
 \setbox\z@\vbox{\hbox{$($}\kern\z@}\bigsize@=2\ht\z@
 \fi
 \normalbaselines\rm\ex@.2326ex\jot5.2\ex@\the\seventeenpoint@}

\font@\rrrrrm=cmr10 scaled\magstep4
\font@\bigtitlefont=cmbx10 scaled\magstep4

\parindent1pc
\normallineskiplimit\p@
\newdimen\indenti \indenti=2pc
\def\pageheight#1{\vsize#1}
\def\pagewidth#1{\hsize#1%
   \captionwidth@\hsize \advance\captionwidth@-2\indenti}
\pagewidth{30pc} \pageheight{47pc}
\def\topmatter{%
 \ifx\undefined\msafam
 \else\font@\eightmsa=msam8 \font@\sixmsa=msam6
   \ifsyntax@\else \addto\tenpoint{\textfont\msafam=\tenmsa
              \scriptfont\msafam=\sevenmsa \scriptscriptfont\msafam=\fivemsa}%
     \addto\eightpoint{\textfont\msafam=\eightmsa \scriptfont\msafam=\sixmsa
              \scriptscriptfont\msafam=\fivemsa}%
   \fi
 \fi
 \ifx\undefined\msbfam
 \else\font@\eightmsb=msbm8 \font@\sixmsb=msbm6
   \ifsyntax@\else \addto\tenpoint{\textfont\msbfam=\tenmsb
         \scriptfont\msbfam=\sevenmsb \scriptscriptfont\msbfam=\fivemsb}%
     \addto\eightpoint{\textfont\msbfam=\eightmsb \scriptfont\msbfam=\sixmsb
         \scriptscriptfont\msbfam=\fivemsb}%
   \fi
 \fi
 \ifx\undefined\eufmfam
 \else \font@\eighteufm=eufm8 \font@\sixeufm=eufm6
   \ifsyntax@\else \addto\tenpoint{\textfont\eufmfam=\teneufm
       \scriptfont\eufmfam=\seveneufm \scriptscriptfont\eufmfam=\fiveeufm}%
     \addto\eightpoint{\textfont\eufmfam=\eighteufm
       \scriptfont\eufmfam=\sixeufm \scriptscriptfont\eufmfam=\fiveeufm}%
   \fi
 \fi
 \ifx\undefined\eufbfam
 \else \font@\eighteufb=eufb8 \font@\sixeufb=eufb6
   \ifsyntax@\else \addto\tenpoint{\textfont\eufbfam=\teneufb
      \scriptfont\eufbfam=\seveneufb \scriptscriptfont\eufbfam=\fiveeufb}%
    \addto\eightpoint{\textfont\eufbfam=\eighteufb
      \scriptfont\eufbfam=\sixeufb \scriptscriptfont\eufbfam=\fiveeufb}%
   \fi
 \fi
 \ifx\undefined\eusmfam
 \else \font@\eighteusm=eusm8 \font@\sixeusm=eusm6
   \ifsyntax@\else \addto\tenpoint{\textfont\eusmfam=\teneusm
       \scriptfont\eusmfam=\seveneusm \scriptscriptfont\eusmfam=\fiveeusm}%
     \addto\eightpoint{\textfont\eusmfam=\eighteusm
       \scriptfont\eusmfam=\sixeusm \scriptscriptfont\eusmfam=\fiveeusm}%
   \fi
 \fi
 \ifx\undefined\eusbfam
 \else \font@\eighteusb=eusb8 \font@\sixeusb=eusb6
   \ifsyntax@\else \addto\tenpoint{\textfont\eusbfam=\teneusb
       \scriptfont\eusbfam=\seveneusb \scriptscriptfont\eusbfam=\fiveeusb}%
     \addto\eightpoint{\textfont\eusbfam=\eighteusb
       \scriptfont\eusbfam=\sixeusb \scriptscriptfont\eusbfam=\fiveeusb}%
   \fi
 \fi
 \ifx\undefined\eurmfam
 \else \font@\eighteurm=eurm8 \font@\sixeurm=eurm6
   \ifsyntax@\else \addto\tenpoint{\textfont\eurmfam=\teneurm
       \scriptfont\eurmfam=\seveneurm \scriptscriptfont\eurmfam=\fiveeurm}%
     \addto\eightpoint{\textfont\eurmfam=\eighteurm
       \scriptfont\eurmfam=\sixeurm \scriptscriptfont\eurmfam=\fiveeurm}%
   \fi
 \fi
 \ifx\undefined\eurbfam
 \else \font@\eighteurb=eurb8 \font@\sixeurb=eurb6
   \ifsyntax@\else \addto\tenpoint{\textfont\eurbfam=\teneurb
       \scriptfont\eurbfam=\seveneurb \scriptscriptfont\eurbfam=\fiveeurb}%
    \addto\eightpoint{\textfont\eurbfam=\eighteurb
       \scriptfont\eurbfam=\sixeurb \scriptscriptfont\eurbfam=\fiveeurb}%
   \fi
 \fi
 \ifx\undefined\cmmibfam
 \else \font@\eightcmmib=cmmib8 \font@\sixcmmib=cmmib6
   \ifsyntax@\else \addto\tenpoint{\textfont\cmmibfam=\tencmmib
       \scriptfont\cmmibfam=\sevencmmib \scriptscriptfont\cmmibfam=\fivecmmib}%
    \addto\eightpoint{\textfont\cmmibfam=\eightcmmib
       \scriptfont\cmmibfam=\sixcmmib \scriptscriptfont\cmmibfam=\fivecmmib}%
   \fi
 \fi
 \ifx\undefined\cmbsyfam
 \else \font@\eightcmbsy=cmbsy8 \font@\sixcmbsy=cmbsy6
   \ifsyntax@\else \addto\tenpoint{\textfont\cmbsyfam=\tencmbsy
      \scriptfont\cmbsyfam=\sevencmbsy \scriptscriptfont\cmbsyfam=\fivecmbsy}%
    \addto\eightpoint{\textfont\cmbsyfam=\eightcmbsy
      \scriptfont\cmbsyfam=\sixcmbsy \scriptscriptfont\cmbsyfam=\fivecmbsy}%
   \fi
 \fi
 \let\topmatter\relax}
\def\chapterno@{\uppercase\expandafter{\romannumeral\chaptercount@}}
\newcount\chaptercount@
\def\chapter{\nofrills@{\afterassignment\chapterno@
                        CHAPTER \global\chaptercount@=}\chapter@
 \DNii@##1{\leavevmode\hskip-\leftskip
   \rlap{\vbox to\z@{\vss\centerline{\eightpoint
   \chapter@##1\unskip}\baselineskip2pc\null}}\hskip\leftskip
   \nofrills@false}%
 \FN@\next@}
\newbox\titlebox@

\def\title{\nofrills@{\relax}\title@%
 \DNii@##1\endtitle{\global\setbox\titlebox@\vtop{\tenpoint\bf
 \raggedcenter@\ignorespaces
 \baselineskip1.3\baselineskip\title@{##1}\endgraf}%
 \ifmonograph@ \edef\next{\the\leftheadtoks}\ifx\next\empty
    \leftheadtext{##1}\fi
 \fi
 \edef\next{\the\rightheadtoks}\ifx\next\empty \rightheadtext{##1}\fi
 }\FN@\next@}
\newbox\authorbox@
\def\author#1\endauthor{\global\setbox\authorbox@
 \vbox{\tenpoint\smc\raggedcenter@\ignorespaces
 #1\endgraf}\relaxnext@ \edef\next{\the\leftheadtoks}%
 \ifx\next\empty\leftheadtext{#1}\fi}
\newbox\affilbox@
\def\affil#1\endaffil{\global\setbox\affilbox@
 \vbox{\tenpoint\raggedcenter@\ignorespaces#1\endgraf}}
\newcount\addresscount@
\addresscount@\z@
\def\address#1\endaddress{\global\advance\addresscount@\@ne
  \expandafter\gdef\csname address\number\addresscount@\endcsname
  {\vskip12\p@ minus6\p@\noindent\eightpoint\smc\ignorespaces#1\par}}
\def\email{\nofrills@{\eightpoint{\it E-mail\/}:\enspace}\email@
  \DNii@##1\endemail{%
  \expandafter\gdef\csname email\number\addresscount@\endcsname
  {\def\usualspace{{\it\enspace}}\smallskip\noindent\eightpoint\email@
  \ignorespaces##1\par}}%
 \FN@\next@}
\def\thedate@{}
\def\date#1\enddate{\gdef\thedate@{\tenpoint\ignorespaces#1\unskip}}
\def\thethanks@{}
\def\thanks#1\endthanks{\gdef\thethanks@{\eightpoint\ignorespaces#1.\unskip}}
\def\thekeywords@{}
\def\keywords{\nofrills@{{\it Key words and phrases.\enspace}}\keywords@
 \DNii@##1\endkeywords{\def\thekeywords@{\def\usualspace{{\it\enspace}}%
 \eightpoint\keywords@\ignorespaces##1\unskip.}}%
 \FN@\next@}
\def\thesubjclass@{}
\def\subjclass{\nofrills@{{\rm2020 {\it Mathematics Subject
   Classification\/}.\enspace}}\subjclass@
 \DNii@##1\endsubjclass{\def\thesubjclass@{\def\usualspace
  {{\rm\enspace}}\eightpoint\subjclass@\ignorespaces##1\unskip.}}%
 \FN@\next@}
\newbox\abstractbox@
\def\abstract{\nofrills@{{\smc Abstract.\enspace}}\abstract@
 \DNii@{\setbox\abstractbox@\vbox\bgroup\noindent$$\vbox\bgroup
  \def\envir@{abstract}\advance\hsize-2\indenti
  \usualspace@{{\enspace}}\eightpoint \noindent\abstract@\ignorespaces}%
 \FN@\next@}
\def\endabstract{\par\unskip\egroup$$\egroup}
\def\widestnumber#1#2{\begingroup\let\head\null\let\subhead\empty
   \let\subsubhead\subhead
   \ifx#1\head\global\setbox\tocheadbox@\hbox{#2.\enspace}%
   \else\ifx#1\subhead\global\setbox\tocsubheadbox@\hbox{#2.\enspace}%
   \else\ifx#1\key\bgroup\let\endrefitem@\egroup
        \key#2\endrefitem@\global\refindentwd\wd\keybox@
   \else\ifx#1\no\bgroup\let\endrefitem@\egroup
        \no#2\endrefitem@\global\refindentwd\wd\nobox@
   \else\ifx#1\page\global\setbox\pagesbox@\hbox{\quad\bf#2}%
   \else\ifx#1\item\setboxz@h{#2}\global\rosteritemwd\wdz@
        \global\advance\rosteritemwd by.5\parindent
   \else\message{\string\widestnumber is not defined for this option
   (\string#1)}%
\fi\fi\fi\fi\fi\fi\endgroup}
\newif\ifmonograph@
\def\Monograph{\monograph@true \let\headmark\rightheadtext
  \let\varindent@\indent \def\headfont@{\bf}\def\proclaimheadfont@{\smc}%
  \def\demofont@{\smc}}
\let\varindent@\indent

\newbox\tocheadbox@    \newbox\tocsubheadbox@
\newbox\tocbox@
\def\toc{\toc@{Contents}}
\def\newtocdefs{%
   \def \title##1\endtitle
       {\penaltyandskip@\z@\smallskipamount
        \hangindent\wd\tocheadbox@\noindent{\bf##1}}%
   \def \chapter##1{%
        Chapter \uppercase\expandafter{\romannumeral##1.\unskip}\enspace}%
   \def \specialhead##1\endspecialhead
       {\par\hangindent\wd\tocheadbox@ \noindent##1\par}%
   \def \head##1 ##2\endhead
       {\par\hangindent\wd\tocheadbox@ \noindent
        \if\notempty{##1}\hbox to\wd\tocheadbox@{\hfil##1\enspace}\fi
        ##2\par}%
   \def \subhead##1 ##2\endsubhead
       {\par\vskip-\parskip {\normalbaselines
        \advance\leftskip\wd\tocheadbox@
        \hangindent\wd\tocsubheadbox@ \noindent
        \if\notempty{##1}\hbox to\wd\tocsubheadbox@{##1\unskip\hfil}\fi
         ##2\par}}%
   \def \subsubhead##1 ##2\endsubsubhead
       {\par\vskip-\parskip {\normalbaselines
        \advance\leftskip\wd\tocheadbox@
        \hangindent\wd\tocsubheadbox@ \noindent
        \if\notempty{##1}\hbox to\wd\tocsubheadbox@{##1\unskip\hfil}\fi
        ##2\par}}}
\def\toc@#1{\relaxnext@
   \def\page##1%
       {\unskip\penalty0\null\hfil
        \rlap{\hbox to\wd\pagesbox@{\quad\hfil##1}}\hfilneg\penalty\@M}%
 \DN@{\ifx\next\nofrills\DN@\nofrills{\nextii@}%
      \else\DN@{\nextii@{{#1}}}\fi
      \next@}%
 \DNii@##1{%
\ifmonograph@\bgroup\else\setbox\tocbox@\vbox\bgroup
   \centerline{\headfont@\ignorespaces##1\unskip}\nobreak
   \vskip\belowheadskip \fi
   \setbox\tocheadbox@\hbox{0.\enspace}%
   \setbox\tocsubheadbox@\hbox{0.0.\enspace}%
   \leftskip\indenti \rightskip\leftskip
   \setbox\pagesbox@\hbox{\bf\quad000}\advance\rightskip\wd\pagesbox@
   \newtocdefs
 }%
 \FN@\next@}
\def\endtoc{\par\egroup}
\let\pretitle\relax
\let\preauthor\relax
\let\preaffil\relax
\let\predate\relax
\let\preabstract\relax
\let\prepaper\relax
\def\dedicatory #1\enddedicatory{\def\preabstract{{\medskip
  \eightpoint\it \raggedcenter@#1\endgraf}}}
\def\thetranslator@{}
\def\translator#1\endtranslator{\def\thetranslator@{\nobreak\medskip
 \line{\eightpoint\hfil Translated by \uppercase{#1}\qquad\qquad}\nobreak}}
\outer\def\endtopmatter{\runaway@{abstract}%
 \edef\next{\the\leftheadtoks}\ifx\next\empty
  \expandafter\leftheadtext\expandafter{\the\rightheadtoks}\fi
 \ifmonograph@\else
   \ifx\thesubjclass@\empty\else \makefootnote@{}{\thesubjclass@}\fi
   \ifx\thekeywords@\empty\else \makefootnote@{}{\thekeywords@}\fi
   \ifx\thethanks@\empty\else \makefootnote@{}{\thethanks@}\fi
 \fi
  \pretitle
  \ifmonograph@ \topskip7pc \else \topskip4pc \fi
  \box\titlebox@
  \topskip10pt
  \preauthor
  \ifvoid\authorbox@\else \vskip2.5pc plus1pc \unvbox\authorbox@\fi
  \preaffil
  \ifvoid\affilbox@\else \vskip1pc plus.5pc \unvbox\affilbox@\fi
  \predate
  \ifx\thedate@\empty\else \vskip1pc plus.5pc \line{\hfil\thedate@\hfil}\fi
  \preabstract
  \ifvoid\abstractbox@\else \vskip1.5pc plus.5pc \unvbox\abstractbox@ \fi
  \ifvoid\tocbox@\else\vskip1.5pc plus.5pc \unvbox\tocbox@\fi
  \prepaper
  \vskip2pc plus1pc
}
\def\document{\let\fontlist@\relax\let\alloclist@\relax
  \tenpoint}

\newskip\aboveheadskip       \aboveheadskip1.8\bigskipamount
\newdimen\belowheadskip      \belowheadskip1.8\medskipamount

\def\headfont@{\smc}
\def\penaltyandskip@#1#2{\relax\ifdim\lastskip<#2\relax\removelastskip
      \ifnum#1=\z@\else\penalty@#1\relax\fi\vskip#2%
  \else\ifnum#1=\z@\else\penalty@#1\relax\fi\fi}
\def\nobreak{\penalty\@M
  \ifvmode\def\penalty@{\let\penalty@\penalty\count@@@}%
  \everypar{\let\penalty@\penalty\everypar{}}\fi}
\let\penalty@\penalty
\def\heading#1\endheading{\head#1\endhead}

\def\specialheadfont@{\bf}
\outer\def\specialhead{\par\penaltyandskip@{-200}\aboveheadskip
  \begingroup\interlinepenalty\@M\rightskip\z@ plus\hsize \let\\\linebreak
  \specialheadfont@\noindent\ignorespaces}
\def\endspecialhead{\par\endgroup\nobreak\vskip\belowheadskip}
\let\headmark\eat@
\newskip\subheadskip       \subheadskip\medskipamount
\def\subheadfont@{\bf}
\outer\def\subhead{\nofrills@{.\enspace}\subhead@
 \DNii@##1\endsubhead{\par\penaltyandskip@{-100}\subheadskip
  \varindent@{\usualspace@{{\subheadfont@\enspace}}%
 \subheadfont@\ignorespaces##1\unskip\subhead@}\ignorespaces}%
 \FN@\next@}
\outer\def\subsubhead{\nofrills@{.\enspace}\subsubhead@
 \DNii@##1\endsubsubhead{\par\penaltyandskip@{-50}\medskipamount
      {\usualspace@{{\it\enspace}}%
  \it\ignorespaces##1\unskip\subsubhead@}\ignorespaces}%
 \FN@\next@}
\def\proclaimheadfont@{\bf}
\outer\def\proclaim{\runaway@{proclaim}\def\envir@{proclaim}%
  \nofrills@{.\enspace}\proclaim@
 \DNii@##1{\penaltyandskip@{-100}\medskipamount\varindent@
   \usualspace@{{\proclaimheadfont@\enspace}}\proclaimheadfont@
   \ignorespaces##1\unskip\proclaim@
  \sl\ignorespaces}%
 \FN@\next@}
\outer\def\endproclaim{\let\envir@\relax\par\rm
  \penaltyandskip@{55}\medskipamount}
\def\demoheadfont@{\it}
\def\demo{\runaway@{proclaim}\nofrills@{.\enspace}\demo@
     \DNii@##1{\par\penaltyandskip@\z@\medskipamount
  {\usualspace@{{\demoheadfont@\enspace}}%
  \varindent@\demoheadfont@\ignorespaces##1\unskip\demo@}\rm
  \ignorespaces}\FN@\next@}
\def\enddemo{\par\medskip}
\def\qed{\ifhmode\unskip\nobreak\fi\quad\ifmmode\square\else$\m@th\square$\fi}
\let\remark\demo
\let\endremark\enddemo
\def\definition{\runaway@{proclaim}%
  \nofrills@{.\demoheadfont@\enspace}\definition@
        \DNii@##1{\penaltyandskip@{-100}\medskipamount
        {\usualspace@{{\demoheadfont@\enspace}}%
        \varindent@\demoheadfont@\ignorespaces##1\unskip\definition@}%
        \rm \ignorespaces}\FN@\next@}


\newdimen\rosteritemwd
\newcount\rostercount@
\newif\iffirstitem@
\let\plainitem@\item
\newtoks\everypartoks@
\def\par@{\everypartoks@\expandafter{\the\everypar}\everypar{}}
\def\roster{\edef\leftskip@{\leftskip\the\leftskip}%
 \relaxnext@
 \rostercount@\z@  
 \def\item{\FN@\rosteritem@}%
 \DN@{\ifx\next\runinitem\let\next@\nextii@\else
  \let\next@\nextiii@\fi\next@}%
 \DNii@\runinitem  
  {\unskip  
   \DN@{\ifx\next[\let\next@\nextii@\else
    \ifx\next"\let\next@\nextiii@\else\let\next@\nextiv@\fi\fi\next@}%
   \DNii@[####1]{\rostercount@####1\relax
    \enspace{\rm(\number\rostercount@)}~\ignorespaces}%
   \def\nextiii@"####1"{\enspace{\rm####1}~\ignorespaces}%
   \def\nextiv@{\enspace{\rm(1)}\rostercount@\@ne~}%
   \par@\firstitem@false  
   \FN@\next@}%
 \def\nextiii@{\par\par@  
  \penalty\@m\smallskip\vskip-\parskip
  \firstitem@true}%
 \FN@\next@}
\def\rosteritem@{\iffirstitem@\firstitem@false\else\par\vskip-\parskip\fi
 \leftskip3\parindent\noindent  
 \DNii@[##1]{\rostercount@##1\relax
  \llap{\hbox to2.5\parindent{\hss\rm(\number\rostercount@)}%
   \hskip.5\parindent}\ignorespaces}%
 \def\nextiii@"##1"{%
  \llap{\hbox to2.5\parindent{\hss\rm##1}\hskip.5\parindent}\ignorespaces}%
 \def\nextiv@{\advance\rostercount@\@ne
  \llap{\hbox to2.5\parindent{\hss\rm(\number\rostercount@)}%
   \hskip.5\parindent}}%
 \ifx\next[\let\next@\nextii@\else\ifx\next"\let\next@\nextiii@\else
  \let\next@\nextiv@\fi\fi\next@}

\newif\ifnextRunin@
\def\endroster{\relaxnext@
 \par\leftskip@  
 \penalty-50 \vskip-\parskip\smallskip  
 \DN@{\ifx\next\Runinitem\let\next@\relax
  \else\nextRunin@false\let\item\plainitem@  
   \ifx\next\par 
    \DN@\par{\everypar\expandafter{\the\everypartoks@}}%
   \else  
    \DN@{\noindent\everypar\expandafter{\the\everypartoks@}}%
  \fi\fi\next@}%
 \FN@\next@}
\newcount\rosterhangafter@
\def\Runinitem#1\roster\runinitem{\relaxnext@
 \rostercount@\z@ 
 \def\item{\FN@\rosteritem@}%
 \def\runinitem@{#1}%
 \DN@{\ifx\next[\let\next\nextii@\else\ifx\next"\let\next\nextiii@
  \else\let\next\nextiv@\fi\fi\next}%
 \DNii@[##1]{\rostercount@##1\relax
  \def\item@{{\rm(\number\rostercount@)}}\nextv@}%
 \def\nextiii@"##1"{\def\item@{{\rm##1}}\nextv@}%
 \def\nextiv@{\advance\rostercount@\@ne
  \def\item@{{\rm(\number\rostercount@)}}\nextv@}%
 \def\nextv@{\setbox\z@\vbox  
  {\ifnextRunin@\noindent\fi  
  \runinitem@\unskip\enspace\item@~\par  
  \global\rosterhangafter@\prevgraf}%
  \firstitem@false  
  \ifnextRunin@\else\par\fi  
  \hangafter\rosterhangafter@\hangindent3\parindent
  \ifnextRunin@\noindent\fi  
  \runinitem@\unskip\enspace 
  \item@~\ifnextRunin@\else\par@\fi  
  \nextRunin@true\ignorespaces}%
 \FN@\next@}
\def\footmarkform@#1{$\m@th^{#1}$}
\let\thefootnotemark\footmarkform@
\def\makefootnote@#1#2{\insert\footins
 {\interlinepenalty\interfootnotelinepenalty
 \eightpoint\splittopskip\ht\strutbox\splitmaxdepth\dp\strutbox
 \floatingpenalty\@MM\leftskip\z@\rightskip\z@\spaceskip\z@\xspaceskip\z@
 \leavevmode{#1}\footstrut\ignorespaces#2\unskip\lower\dp\strutbox
 \vbox to\dp\strutbox{}}}
\newcount\footmarkcount@
\footmarkcount@\z@
\def\footnotemark{\let\@sf\empty\relaxnext@
 \ifhmode\edef\@sf{\spacefactor\the\spacefactor}\/\fi
 \DN@{\ifx[\next\let\next@\nextii@\else
  \ifx"\next\let\next@\nextiii@\else
  \let\next@\nextiv@\fi\fi\next@}%
 \DNii@[##1]{\footmarkform@{##1}\@sf}%
 \def\nextiii@"##1"{{##1}\@sf}%
 \def\nextiv@{\iffirstchoice@\global\advance\footmarkcount@\@ne\fi
  \footmarkform@{\number\footmarkcount@}\@sf}%
 \FN@\next@}
\def\footnotetext{\relaxnext@
 \DN@{\ifx[\next\let\next@\nextii@\else
  \ifx"\next\let\next@\nextiii@\else
  \let\next@\nextiv@\fi\fi\next@}%
 \DNii@[##1]##2{\makefootnote@{\footmarkform@{##1}}{##2}}%
 \def\nextiii@"##1"##2{\makefootnote@{##1}{##2}}%
 \def\nextiv@##1{\makefootnote@{\footmarkform@{\number\footmarkcount@}}{##1}}%
 \FN@\next@}
\def\footnote{\let\@sf\empty\relaxnext@
 \ifhmode\edef\@sf{\spacefactor\the\spacefactor}\/\fi
 \DN@{\ifx[\next\let\next@\nextii@\else
  \ifx"\next\let\next@\nextiii@\else
  \let\next@\nextiv@\fi\fi\next@}%
 \DNii@[##1]##2{\footnotemark[##1]\footnotetext[##1]{##2}}%
 \def\nextiii@"##1"##2{\footnotemark"##1"\footnotetext"##1"{##2}}%
 \def\nextiv@##1{\footnotemark\footnotetext{##1}}%
 \FN@\next@}
\def\adjustfootnotemark#1{\advance\footmarkcount@#1\relax}
\def\footnoterule{\kern-3\p@
  \hrule width 5pc\kern 2.6\p@} 
\def\captionfont@{\smc}
\def\topcaption#1#2\endcaption{%
  {\dimen@\hsize \advance\dimen@-\captionwidth@
   \rm\raggedcenter@ \advance\leftskip.5\dimen@ \rightskip\leftskip
  {\captionfont@#1}%
  \if\notempty{#2}.\enspace\ignorespaces#2\fi
  \endgraf}\nobreak\bigskip}
\def\botcaption#1#2\endcaption{%
  \nobreak\bigskip
  \setboxz@h{\captionfont@#1\if\notempty{#2}.\enspace\rm#2\fi}%
  {\dimen@\hsize \advance\dimen@-\captionwidth@
   \leftskip.5\dimen@ \rightskip\leftskip
   \noindent \ifdim\wdz@>\captionwidth@ 
   \else\hfil\fi 
  {\captionfont@#1}\if\notempty{#2}.\enspace\rm#2\fi\endgraf}}
\def\@ins{\par\begingroup\def\vspace##1{\vskip##1\relax}%
  \def\captionwidth##1{\captionwidth@##1\relax}%
  \setbox\z@\vbox\bgroup} 
\def\block{\RIfMIfI@\nondmatherr@\block\fi
       \else\ifvmode\vskip\abovedisplayskip\noindent\fi
        $$\def\endblock{\par\egroup$$}\fi
  \vbox\bgroup\advance\hsize-2\indenti\noindent}
\def\endblock{\par\egroup}
\def\cite#1{{\rm[{\citefont@\m@th#1}]}}
\def\citefont@{\rm}
\def\refsfont@{\eightpoint}
\outer\def\Refs{\runaway@{proclaim}%
 \relaxnext@ \DN@{\ifx\next\nofrills\DN@\nofrills{\nextii@}\else
  \DN@{\nextii@{References}}\fi\next@}%
 \DNii@##1{\penaltyandskip@{-200}\aboveheadskip
  \line{\hfil\headfont@\ignorespaces##1\unskip\hfil}\nobreak
  \vskip\belowheadskip
  \begingroup\refsfont@\sfcode`.=\@m}%
 \FN@\next@}
\def\endRefs{\par\endgroup}
\newbox\nobox@            \newbox\keybox@           \newbox\bybox@
\newbox\paperbox@         \newbox\paperinfobox@     \newbox\jourbox@
\newbox\volbox@           \newbox\issuebox@         \newbox\yrbox@
\newbox\pagesbox@         \newbox\bookbox@          \newbox\bookinfobox@
\newbox\publbox@          \newbox\publaddrbox@      \newbox\finalinfobox@
\newbox\edsbox@           \newbox\langbox@
\newif\iffirstref@        \newif\iflastref@
\newif\ifprevjour@        \newif\ifbook@            \newif\ifprevinbook@
\newif\ifquotes@          \newif\ifbookquotes@      \newif\ifpaperquotes@
\newdimen\bysamerulewd@
\setboxz@h{\refsfont@\kern3em}
\bysamerulewd@\wdz@
\newdimen\refindentwd
\setboxz@h{\refsfont@ 00. }
\refindentwd\wdz@
\outer\def\ref{\begingroup \noindent\hangindent\refindentwd
 \firstref@true \def\nofrills{\def\refkern@{\kern3sp}}%
 \ref@}
\def\ref@{\book@false \bgroup\let\endrefitem@\egroup \ignorespaces}
\def\moreref{\endrefitem@\endref@\firstref@false\ref@}%
\def\transl{\endrefitem@\endref@\firstref@false
  \book@false
  \prepunct@
  \setboxz@h\bgroup \aftergroup\unhbox\aftergroup\z@
    \def\endrefitem@{\unskip\refkern@\egroup}\ignorespaces}%
\def\emptyifempty@{\dimen@\wd\currbox@
  \advance\dimen@-\wd\z@ \advance\dimen@-.1\p@
  \ifdim\dimen@<\z@ \setbox\currbox@\copy\voidb@x \fi}
\let\refkern@\relax
\def\endrefitem@{\unskip\refkern@\egroup
  \setboxz@h{\refkern@}\emptyifempty@}\ignorespaces
\def\refdef@#1#2#3{\edef\next@{\noexpand\endrefitem@
  \let\noexpand\currbox@\csname\expandafter\eat@\string#1box@\endcsname
    \noexpand\setbox\noexpand\currbox@\hbox\bgroup}%
  \toks@\expandafter{\next@}%
  \if\notempty{#2#3}\toks@\expandafter{\the\toks@
  \def\endrefitem@{\unskip#3\refkern@\egroup
  \setboxz@h{#2#3\refkern@}\emptyifempty@}#2}\fi
  \toks@\expandafter{\the\toks@\ignorespaces}%
  \edef#1{\the\toks@}}
\refdef@\no{}{. }
\refdef@\key{[\m@th}{] }
\refdef@\by{}{}
\def\bysame{\by\hbox to\bysamerulewd@{\hrulefill}\thinspace
   \kern0sp}
\def\manyby{\message{\string\manyby is no longer necessary; \string\by
  can be used instead, starting with version 2.0 of \styname.STY}\by}
\refdef@\paper{\ifpaperquotes@``\fi\it}{}
\refdef@\paperinfo{}{}
\def\jour{\endrefitem@\let\currbox@\jourbox@
  \setbox\currbox@\hbox\bgroup
  \def\endrefitem@{\unskip\refkern@\egroup
    \setboxz@h{\refkern@}\emptyifempty@
    \ifvoid\jourbox@\else\prevjour@true\fi}%
\ignorespaces}
\refdef@\vol{\ifbook@\else\bf\fi}{}
\refdef@\issue{no. }{}
\refdef@\yr{}{}
\refdef@\pages{}{}
\def\page{\endrefitem@\def\pp@{\def\pp@{pp.~}p.~}\let\currbox@\pagesbox@
  \setbox\currbox@\hbox\bgroup\ignorespaces}
\def\pp@{pp.~}
\def\book{\endrefitem@ \let\currbox@\bookbox@
 \setbox\currbox@\hbox\bgroup\def\endrefitem@{\unskip\refkern@\egroup
  \setboxz@h{\ifbookquotes@``\fi}\emptyifempty@
  \ifvoid\bookbox@\else\book@true\fi}%
  \ifbookquotes@``\fi\it\ignorespaces}
\def\inbook{\endrefitem@
  \let\currbox@\bookbox@\setbox\currbox@\hbox\bgroup
  \def\endrefitem@{\unskip\refkern@\egroup
  \setboxz@h{\ifbookquotes@``\fi}\emptyifempty@
  \ifvoid\bookbox@\else\book@true\previnbook@true\fi}%
  \ifbookquotes@``\fi\ignorespaces}
\refdef@\eds{(}{, eds.)}
\def\ed{\endrefitem@\let\currbox@\edsbox@
 \setbox\currbox@\hbox\bgroup
 \def\endrefitem@{\unskip, ed.)\refkern@\egroup
  \setboxz@h{(, ed.)}\emptyifempty@}(\ignorespaces}
\refdef@\bookinfo{}{}
\refdef@\publ{}{}
\refdef@\publaddr{}{}
\refdef@\finalinfo{}{}
\refdef@\lang{(}{)}

\let\refdef@\relax 
\def\ppunbox@#1{\ifvoid#1\else\prepunct@\unhbox#1\fi}
\def\nocomma@#1{\ifvoid#1\else\changepunct@3\prepunct@\unhbox#1\fi}
\def\changepunct@#1{\ifnum\lastkern<3 \unkern\kern#1sp\fi}
\def\prepunct@{\count@\lastkern\unkern
  \ifnum\lastpenalty=0
    \let\penalty@\relax
  \else
    \edef\penalty@{\penalty\the\lastpenalty\relax}%
  \fi
  \unpenalty
  \let\refspace@\ \ifcase\count@,
\or;\or.\or 
  \or\let\refspace@\relax
  \else,\fi
  \ifquotes@''\quotes@false\fi \penalty@ \refspace@
}
\def\transferpenalty@#1{\dimen@\lastkern\unkern
  \ifnum\lastpenalty=0\unpenalty\let\penalty@\relax
  \else\edef\penalty@{\penalty\the\lastpenalty\relax}\unpenalty\fi
  #1\penalty@\kern\dimen@}
\def\endref{\endrefitem@\lastref@true\endref@
  \par\endgroup \prevjour@false \previnbook@false }
\def\endref@{%
\iffirstref@
  \ifvoid\nobox@\ifvoid\keybox@\indent\fi
  \else\hbox to\refindentwd{\hss\unhbox\nobox@}\fi
  \ifvoid\keybox@
  \else\ifdim\wd\keybox@>\refindentwd
         \box\keybox@
       \else\hbox to\refindentwd{\unhbox\keybox@\hfil}\fi\fi
  \kern4sp\ppunbox@\bybox@
\fi 
  \ifvoid\paperbox@
  \else\prepunct@\unhbox\paperbox@
    \ifpaperquotes@\quotes@true\fi\fi
  \ppunbox@\paperinfobox@
  \ifvoid\jourbox@
    \ifprevjour@ \nocomma@\volbox@
      \nocomma@\issuebox@
      \ifvoid\yrbox@\else\changepunct@3\prepunct@(\unhbox\yrbox@
        \transferpenalty@)\fi
      \ppunbox@\pagesbox@
    \fi 
  \else \prepunct@\unhbox\jourbox@
    \nocomma@\volbox@
    \nocomma@\issuebox@
    \ifvoid\yrbox@\else\changepunct@3\prepunct@(\unhbox\yrbox@
      \transferpenalty@)\fi
    \ppunbox@\pagesbox@
  \fi 
  \ifbook@\prepunct@\unhbox\bookbox@ \ifbookquotes@\quotes@true\fi \fi
  \nocomma@\edsbox@
  \ppunbox@\bookinfobox@
  \ifbook@\ifvoid\volbox@\else\prepunct@ vol.~\unhbox\volbox@
  \fi\fi
  \ppunbox@\publbox@ \ppunbox@\publaddrbox@
  \ifbook@ \ppunbox@\yrbox@
    \ifvoid\pagesbox@
    \else\prepunct@\pp@\unhbox\pagesbox@\fi
  \else
    \ifprevinbook@ \ppunbox@\yrbox@
      \ifvoid\pagesbox@\else\prepunct@\pp@\unhbox\pagesbox@\fi
    \fi \fi
  \ppunbox@\finalinfobox@
  \iflastref@
    \ifvoid\langbox@.\ifquotes@''\fi
    \else\changepunct@2\prepunct@\unhbox\langbox@\fi
  \else
    \ifvoid\langbox@\changepunct@1%
    \else\changepunct@3\prepunct@\unhbox\langbox@
      \changepunct@1\fi
  \fi
}
\outer\def\enddocument{%
 \runaway@{proclaim}%
\ifmonograph@ 
\else
 \nobreak
 \thetranslator@
 \count@\z@ \loop\ifnum\count@<\addresscount@\advance\count@\@ne
 \csname address\number\count@\endcsname
 \csname email\number\count@\endcsname
 \repeat
\fi
 \vfill\supereject\end}

\def\headfont@{\headfonts}
\def\proclaimheadfont@{\bf}
\def\specialheadfont@{\bf}
\def\subheadfont@{\bf}
\def\demoheadfont@{\smc}

\newif\ifThisToToc \ThisToTocfalse
\newif\iftocloaded \tocloadedfalse

\def\C@L{\noexpand\Cal}\def\B@B{\noexpand\Bbb}\def\fR@K{\noexpand\frak}
\def\S@{\noexpand\S}\def\P@P{\noexpand\"}
\def\xpar{\\}

\def\writetoc#1{\iftocloaded\ifThisToToc\begingroup\def\totoc{}
  \def\Cal{\noexpand\C@L}\def\Bbb{\noexpand\B@B}
  \def\frak{\noexpand\fR@K}\def\goth{\frak}\def\S{\noexpand\S@}
  \def\"{\noexpand\P@P}
  \def\xpar{\par\penalty100000 }\def\idx##1{##1}\def\\{\xpar}
  \edef\next@{\write\toc{\noindent#1\leaderfill\noexpand\folio\par}}%
  \next@\endgroup\global\ThisToTocfalse\fi\fi}
\def\leaderfill{\leaders\hbox to 1em{\hss.\hss}\hfill}

\newif\ifindexloaded \indexloadedfalse
\def\idx#1{\ifindexloaded\begingroup\def\ign{}\def\it{}\def\/{}%
 \def\smc{}\def\bf{}\def\tt{}%
 \def\Cal{\noexpand\C@L}\def\Bbb{\noexpand\B@B}%
 \def\frak{\noexpand\fR@K}\def\goth{\frak}\def\S{\noexpand\S@}%
  \def\"{\noexpand\P@P}%
 {\edef\next@{\write\index{#1, \noexpand\folio}}\next@}%
 \endgroup\fi{#1}}
\def\ign#1{}

\def\totoc{\global\ThisToToctrue}

\outer\def\head#1\endhead{\par\penaltyandskip@{-200}\aboveheadskip
 {\headfont@\raggedcenter@\interlinepenalty\@M
 \ignorespaces#1\endgraf}\nobreak
 \vskip\belowheadskip
 \headmark{#1}\writetoc{#1}}

\outer\def\chaphead#1\endchaphead{\par\penaltyandskip@{-200}\aboveheadskip
 {\chapheadfonts\raggedcenter@\interlinepenalty\@M
 \ignorespaces#1\endgraf}\nobreak
 \vskip3\belowheadskip
 \headmark{#1}\writetoc{#1}}

\def\folio{{\foliofont@\ifnum\pageno<\z@ \romannumeral-\pageno
 \else\number\pageno \fi}}
\newtoks\leftheadtoks
\newtoks\rightheadtoks

\def\leftheadtext{\nofrills@{\relax}\lht@
  \DNii@##1{\leftheadtoks\expandafter{\lht@{##1}}%
    \mark{\the\leftheadtoks\noexpand\else\the\rightheadtoks}
    \ifsyntax@\setboxz@h{\def\\{\unskip\space\ignorespaces}%
        \headlinefont@##1}\fi}%
  \FN@\next@}
\def\rightheadtext{\nofrills@{\relax}\rht@
  \DNii@##1{\rightheadtoks\expandafter{\rht@{##1}}%
    \mark{\the\leftheadtoks\noexpand\else\the\rightheadtoks}%
    \ifsyntax@\setboxz@h{\def\\{\unskip\space\ignorespaces}%
        \headlinefont@##1}\fi}%
  \FN@\next@}
\def\NoRunningHeads{\global\runheads@false\global\let\headmark\eat@}

\newif\iffirstpage@     \firstpage@true
\newif\ifrunheads@      \runheads@true

\newdimen\fullhsize \fullhsize=\hsize
\newdimen\fullvsize \fullvsize=\vsize
\def\fullline{\hbox to\fullhsize}

\def\pagenumbers{\gdef\folio{\folio@}}

\let\norunningheads\NoRunningHeads
\def\userunningheads{\global\runheads@true}
\norunningheads

\headline={\def\chapter#1{\chapterno@. }%
  \def\\{\unskip\space\ignorespaces}\ifrunheads@\headlinefont@
    \ifodd\pageno\rightheadline \else\leftheadline\fi
   \else\hfil\fi\ifNoRunHeadline\global\NoRunHeadlinefalse\fi}
\let\folio@\folio
\def\foliofont@{\foliofont}
\def\foliofont{\eightrm}
\def\headlinefont@{\headlinefont}
\def\headlinefont{\eightpoint\smc}
\def\leftheadline{\rlap{\folio}\hfill
   \ifNoRunHeadline\else\iftrue\topmark\fi\fi \hfill}
\def\rightheadline{\hfill\ifNoRunHeadline
   \else \expandafter\fi
  \hfill \llap{\folio}}
\footline={{\eightpoint\bottremark}%
   \ifrunheads@\else\hfil{\let\foliofont\tenrm\folio}\fi\hfil}
\def\bottremark{}
 
\newif\ifNoRunHeadline      
\def\norunninghead{\global\NoRunHeadlinetrue}
\norunninghead

\output={\output@}
%
\newif\ifoffset\offsetfalse
\output={\output@}
\def\output@{%
 \ifoffset 
  \ifodd\count0\advance\hoffset by0.5truecm
   \else\advance\hoffset by-0.5truecm\fi\fi
 \shipout\vbox{%
  \makeheadline \pagebody \makefootline }%
 \advancepageno \ifnum\outputpenalty>-\@MM\else\dosupereject\fi}

\def\indexoutput#1{%
  \ifoffset 
   \ifodd\count0\advance\hoffset by0.5truecm
    \else\advance\hoffset by-0.5truecm\fi\fi
  \shipout\vbox{\makeheadline
  \vbox to\fullvsize{\boxmaxdepth\maxdepth%
  \ifvoid\topins\else\unvbox\topins\fi%
  #1 %
  \ifvoid\footins\else 
    \vskip\skip\footins
    \footnoterule
    \unvbox\footins\fi
  \ifr@ggedbottom \kern-\dimen@ \vfil \fi}%
  \baselineskip2pc
  \makefootline}%
 \global\advance\pageno\@ne
 \ifnum\outputpenalty>-\@MM\else\dosupereject\fi}
 
 \newbox\partialpage \newdimen\halfsize \halfsize=0.5\fullhsize
 \advance\halfsize by-0.5em

 \def\begindoublecolumns{\output={\indexoutput{\unvbox255}}%
   \begingroup \def\line{\fullline}
   \output={\global\setbox\partialpage=\vbox{\unvbox255\bigskip}}\eject
   \output={\doublecolumnout}\hsize=\halfsize \vsize=2\fullvsize}
 \def\enddoublecolumns{\output={\balancecolumns}\eject
  \endgroup \pagegoal=\fullvsize%
  \output={\output@}}
\def\doublecolumnout{\splittopskip=\topskip \splitmaxdepth=\maxdepth
  \dimen@=\fullvsize \advance\dimen@ by-\ht\partialpage
  \setbox0=\vsplit255 to \dimen@ \setbox2=\vsplit255 to \dimen@
  \indexoutput{\pagesofar} \unvbox255 \penalty\outputpenalty}
\def\pagesofar{\unvbox\partialpage
  \wd0=\hsize \wd2=\hsize \hbox to\fullhsize{\box0\hfil\box2}}
\def\balancecolumns{\setbox0=\vbox{\unvbox255} \dimen@=\ht0
  \advance\dimen@ by\topskip \advance\dimen@ by-\baselineskip
  \divide\dimen@ by2 \splittopskip=\topskip
  {\vbadness=10000 \loop \global\setbox3=\copy0
    \global\setbox1=\vsplit3 to\dimen@
    \ifdim\ht3>\dimen@ \global\advance\dimen@ by1pt \repeat}
  \setbox0=\vbox to\dimen@{\unvbox1} \setbox2=\vbox to\dimen@{\unvbox3}
  \pagesofar}

\tenpoint
\catcode`\@=\active

\def\smallheadings{\let\chapheadfonts\tenpoint\let\headfonts\tenpoint}

\tenpoint
\catcode`\@=\active

\def\LL{\leavevmode\setbox0=\hbox{L}\hbox to\wd0{\hss\char'40L}}
\def\al{\alpha}

\def\ga{\gamma}


\def\today{\ifcase\month\or
 January\or February\or March\or April\or May\or June\or
 July\or August\or September\or October\or November\or December\fi
 \space\number\day, \number\year}

\def\({\left(}
\def\){\right)}
\def\[{\left[}
\def\]{\right]}

\def\3{\ss}
\catcode`\@=11
\def\dddot#1{\vbox{\ialign{##\crcr
      .\hskip-.5pt.\hskip-.5pt.\crcr\noalign{\kern1.5\p@\nointerlineskip}
      $\hfil\displaystyle{#1}\hfil$\crcr}}}

\newif\iftab@\tab@false
\newif\ifvtab@\vtab@false
\def\tab{\bgroup\tab@true\vtab@false\vst@bfalse\Strich@false%
   \def\\{\global\hline@@false%
     \ifhline@\global\hline@false\global\hline@@true\fi\cr}
   \edef\l@{\the\leftskip}\ialign\bgroup\hskip\l@##\hfil&&##\hfil\cr}
\def\endtab{\cr\egroup\egroup}
\def\vtab{\vtop\bgroup\vst@bfalse\vtab@true\tab@true\Strich@false%
   \bgroup\def\\{\cr}\ialign\bgroup&##\hfil\cr}
\def\endvtab{\cr\egroup\egroup\egroup}
\def\stab{\D@cke0.5pt\null 
 \bgroup\tab@true\vtab@false\vst@bfalse\Strich@true\Let@@\vspace@
 \normalbaselines\offinterlineskip
  \openup\spreadmlines@
 \edef\l@{\the\leftskip}\ialign
 \bgroup\hskip\l@##\hfil&&##\hfil\crcr}
\def\endstab{\crcr\egroup
 \egroup}
\newif\ifvst@b\vst@bfalse
\def\vstab{\D@cke0.5pt\null
 \vtop\bgroup\tab@true\vtab@false\vst@btrue\Strich@true\bgroup\Let@@\vspace@
 \normalbaselines\offinterlineskip
  \openup\spreadmlines@\bgroup}
\def\endvstab{\crcr\egroup\egroup
 \egroup\tab@false\Strich@false}

\newdimen\htstrut@
\htstrut@8.5\p@
\newdimen\htStrut@
\htStrut@12\p@
\newdimen\dpstrut@
\dpstrut@3.5\p@
\newdimen\dpStrut@
\dpStrut@3.5\p@
\def\openup{\afterassignment\@penup\dimen@=}
\def\@penup{\advance\lineskip\dimen@
  \advance\baselineskip\dimen@
  \advance\lineskiplimit\dimen@
  \divide\dimen@ by2
  \advance\htstrut@\dimen@
  \advance\htStrut@\dimen@
  \advance\dpstrut@\dimen@
  \advance\dpStrut@\dimen@}
\def\Let@@{\relax%
    \def\\{\global\hline@@false%
     \ifhline@\global\hline@false\global\hline@@true\fi\cr}%
    \iffalse}\fi}
\def\matrix{\null\,\vcenter\bgroup
 \tab@false\vtab@false\vst@bfalse\Strich@false\Let@@\vspace@
 \normalbaselines\openup\spreadmlines@\ialign
 \bgroup\hfil$\m@th##$\hfil&&\quad\hfil$\m@th##$\hfil\crcr
 \Mathstrut@\crcr\noalign{\kern-\baselineskip}}
\def\endmatrix{\crcr\Mathstrut@\crcr\noalign{\kern-\baselineskip}\egroup
 \egroup\,}
\def\smatrix{\D@cke0.5pt\null\,
 \vcenter\bgroup\tab@false\vtab@false\vst@bfalse\Strich@true\Let@@\vspace@
 \normalbaselines\offinterlineskip
  \openup\spreadmlines@\ialign
 \bgroup\hfil$\m@th##$\hfil&&\quad\hfil$\m@th##$\hfil\crcr}
\def\endsmatrix{\crcr\egroup
 \egroup\,\Strich@false}
\newdimen\D@cke
\def\Dicke#1{\global\D@cke#1}
\newtoks\tabs@\tabs@{&}
\newif\ifStrich@\Strich@false
\newif\iff@rst

\def\Stricherr@{\iftab@\ifvtab@\errmessage{\noexpand\s not allowed
     here. Use \noexpand\vstab!}%
  \else\errmessage{\noexpand\s not allowed here. Use \noexpand\stab!}%
  \fi\else\errmessage{\noexpand\s not allowed
     here. Use \noexpand\smatrix!}\fi}
\def\format{\ifvst@b\else\crcr\fi\egroup\iffalse{\fi\ifnum`}=0 \fi\format@}
\def\format@#1\\{\def\preamble@{#1}%
 \def\Str@chfehlt##1{\ifx##1\s\Stricherr@\fi\ifx##1\\\let\Next\relax%
   \else\let\Next\Str@chfehlt\fi\Next}%
 \def\c{\hfil\noexpand\ifhline@@\hbox{\vrule height\htStrut@%
   depth\dpstrut@ width\z@}\noexpand\fi%
   \ifStrich@\hbox{\vrule height\htstrut@ depth\dpstrut@ width\z@}%
   \fi\iftab@\else$\m@th\fi\the\hashtoks@\iftab@\else$\fi\hfil}%
 \def\r{\hfil\noexpand\ifhline@@\hbox{\vrule height\htStrut@%
   depth\dpstrut@ width\z@}\noexpand\fi%
   \ifStrich@\hbox{\vrule height\htstrut@ depth\dpstrut@ width\z@}%
   \fi\iftab@\else$\m@th\fi\the\hashtoks@\iftab@\else$\fi}%
 \def\l{\noexpand\ifhline@@\hbox{\vrule height\htStrut@%
   depth\dpstrut@ width\z@}\noexpand\fi%
   \ifStrich@\hbox{\vrule height\htstrut@ depth\dpstrut@ width\z@}%
   \fi\iftab@\else$\m@th\fi\the\hashtoks@\iftab@\else$\fi\hfil}%
 \def\s{\ifStrich@\ \the\tabs@\vrule width\D@cke\the\hashtoks@%
          \fi\the\tabs@\ }%
 \def\sa{\ifStrich@\vrule width\D@cke\the\hashtoks@%
            \the\tabs@\ %
            \fi}%
 \def\se{\ifStrich@\ \the\tabs@\vrule width\D@cke\the\hashtoks@\fi}%
 \def\cd{\hfil\noexpand\ifhline@@\hbox{\vrule height\htStrut@%
   depth\dpstrut@ width\z@}\noexpand\fi%
   \ifStrich@\hbox{\vrule height\htstrut@ depth\dpstrut@ width\z@}%
   \fi$\dsize\m@th\the\hashtoks@$\hfil}%
 \def\rd{\hfil\noexpand\ifhline@@\hbox{\vrule height\htStrut@%
   depth\dpstrut@ width\z@}\noexpand\fi%
   \ifStrich@\hbox{\vrule height\htstrut@ depth\dpstrut@ width\z@}%
   \fi$\dsize\m@th\the\hashtoks@$}%
 \def\ld{\noexpand\ifhline@@\hbox{\vrule height\htStrut@%
   depth\dpstrut@ width\z@}\noexpand\fi%
   \ifStrich@\hbox{\vrule height\htstrut@ depth\dpstrut@ width\z@}%
   \fi$\dsize\m@th\the\hashtoks@$\hfil}%
 \ifStrich@\else\Str@chfehlt#1\\\fi%
 \setbox\z@\hbox{\xdef\Preamble@{\preamble@}}\ifnum`{=0 \fi\iffalse}\fi
 \ialign\bgroup\span\Preamble@\crcr}
\newif\ifhline@\hline@false
\newif\ifhline@@\hline@@false
\def\hlinefor#1{\multispan@{\strip@#1 }\leaders\hrule height\D@cke\hfill%
    \global\hline@true\ignorespaces}
\def\Item "#1"{\par\noindent\hangindent2\parindent%
  \hangafter1\setbox0\hbox{\rm#1\enspace}\ifdim\wd0>2\parindent%
  \box0\else\hbox to 2\parindent{\rm#1\hfil}\fi\ignorespaces}
\def\ITEM #1"#2"{\par\noindent\hangafter1\hangindent#1%
  \setbox0\hbox{\rm#2\enspace}\ifdim\wd0>#1%
  \box0\else\hbox to 0pt{\rm#2\hss}\hskip#1\fi\ignorespaces}
\def\item"#1"{\par\noindent\hang%
  \setbox0=\hbox{\rm#1\enspace}\ifdim\wd0>\the\parindent%
  \box0\else\hbox to \parindent{\rm#1\hfil}\enspace\fi\ignorespaces}
\let\plainitem@\item
\catcode`\@=13

\hsize13cm
\vsize19cm
\newdimen\fullhsize
\newdimen\fullvsize
\newdimen\halfsize
\fullhsize13cm
\fullvsize19cm
\halfsize=0.5\fullhsize
\advance\halfsize by-0.5em

\magnification1200

\TagsOnRight

\def\AndrAF{1}
\def\AndrCB{2}
\def\ArmDAB{3}
\def\BiSwAA{4}
\def\HoRuAA{5}
\def\KrMuAD{6}
\def\ReSWAB{7}
\def\StanBI{8}
\def\StatAD{9}
\def\StucAA{10}
\def\WaZuAB{11}
\def\ZudiAA{12}

\def\AA{1}
\def\AS{2}
\def\AL{3}
\def\AF{4}
\def\AM{5}
\def\AN{6}
\def\AO{7}
\def\AQ{8}
\def\AP{9}
\def\AE{10}
\def\AH{11}
\def\AG{12}
\def\AI{13}
\def\AB{14}
\def\AD{15}
\def\AJ{16}
\def\AK{17}
\def\AC{18}
\def\AR{19}

\def\TE{1}
\def\RA{2}
\def\TG{3}
\def\TH{4}
\def\TD{5}
\def\TF{6}
\def\RB{7}
\def\TA{8}
\def\TB{9}
\def\TC{10}

\def\fl#1{\left\lfloor#1\right\rfloor}

\topmatter 
\title A positivity conjecture for a quotient of $q$-binomial coefficients
\endtitle 
\author M.  Gatzweiler and C.~Krattenthaler
\endauthor

\affil 
Fakult\"at f\"ur Mathematik, Universit\"at Wien,\\
Oskar-Morgenstern-Platz~1, A-1090 Vienna, Austria.\\
WWW: \tt http://www.mat.univie.ac.at/\~{}kratt
\endaffil
\address Fakult\"at f\"ur Mathematik, Universit\"at Wien,
Oskar-Morgenstern-Platz~1, A-1090 Vienna, Austria.\newline
http://www.mat.univie.ac.at/\~{}kratt
\endaddress

\subjclass Primary 05A30;
 Secondary 33D05
\endsubjclass
\keywords $q$-binomial coefficients, positivity,
rational $q$-Catalan polynomial
\endkeywords

\dedicatory{Dedicated to the grandmasters of $q$-calculus,
George Andrews and Bruce Berndt}
\enddedicatory

\thanks The first author has been supported by the University of
Vienna in the framework of the Vienna School of Mathematics.
The second author has been partially supported by the Austrian
Science Foundation FWF, grant 10.55776/F1002, in the framework
of the Special Research Program ``Discrete Random Structures:
Enumeration and Scaling Limits"%
\endthanks

\abstract We conjecture that, if the quotient of two $q$-binomial
coefficients with the same top argument
is a polynomial, then it has non-negative coefficients.
We summarise what is known about the conjecture and prove it in two
non-trivial cases. Moreover, we move ahead to extend our conjecture
to D.~Stanton's fake Gaussian sequences.
As a corollary of one of our results we obtain
that a polynomial that is conjectured to
be a cyclic sieving polynomial for Kreweras words
[S.~Hopkins and M. Rubey, 
{\it Selecta Math\. (N.S.)} {\bf28} (2022), Paper No.~10]
is  indeed a polynomial with non-negative integer coefficients.
\endabstract
\endtopmatter
\document

We use the standard $q$-notations $[n]_q:=\frac {1-q^n} {1-q}$,
$$
[n]_q!:=[n]_q\,[n-1]_q\cdots[1]_q, \quad \text{for }n\ge1,
\quad \text{and $[0]_q!:=1$,}
$$
and
$$
\bmatrix n\\k\endbmatrix_q:=\cases 
\frac {[n]_q!} {[k]_q!\,[n-k]_q!},&\text{if }0\le k\le n,\\
0,&\text{otherwise,}
\endcases
$$
to define {\it $q$-factorials} and {\it $q$-binomial coefficients},
respectively.

The purpose of this note is to put forward the following conjecture.

\proclaim{Conjecture \TE}
Let $n,k,l$ be non-negative integers. If
$$
\frac {\bmatrix n\\k\endbmatrix_q} {\bmatrix n\\l\endbmatrix_q}
=\frac {[l]_q!\,[n-l]_q!} {[k]_q!\,[n-k]_q!} 
\tag\AA
$$
is a polynomial\/ {\rm(}in $q${\rm)} then it has non-negative coefficients.
\endproclaim

\remark{Remark \RA}
(1)
It is easy to see that it is enough to consider
the range $1\le l<k\le \frac {n} {2}$.

\medskip
(2) In evidence for this conjecture, we point out that it is
trivially true for $l=k$ and $l=k-1$, well known for $l=0$ and $l=1$,
implicitly known for $l=k-2$, and
that we prove it for $l=k-3$ and $l=2$.
These cases are treated in Cases~1--7 below and form the bulk of
this paper.
Moreover,
we have verified Conjecture~\TE\ for all $n\le 400$ with the help
of {\sl Mathematica}.

\medskip
(3) As Sam Hopkins pointed out to us as a reaction on an earlier version of
this paper, our Conjecture~\TE\ has some overlap with an unpublished
conjecture of Stanton on {\it fake Gaussian sequences}~\cite{\StatAD, Conj.~1}.
It says:

\medskip
{\it If
$$
\prod _{i=1} ^{n}\frac {(1-q^{m+i})^{a_i}}
{(1-q^{i})^{a_i}}
\tag\AS$$
is a polynomial for all positive integers~$m$, where $n$ is
another positive integer,
and\linebreak $(a_1,a_2,\dots,a_n)$ is a symmetric
sequence of non-negative integers, then it has non-negative coefficients
for all~$m$.
}

\medskip
Assuming $k>l$, if in the expression in (\AS)
we replace $n$ by $n-k+l$, choose $m=k-l$ and $(a_1,a_2,\dots,a_n)
=(0^l,1^{n-k-l},0^l)$, where $\al^s$ stands for $\al$ repeated $s$~times,
then we obtain~(\AA).\footnote{As we shall see in Case~5 below, for
$l=0$ the expression in~(\AA) reduces to a $q$-binomial coefficient,
which is also known as {\it Gau\ss ian polynomial\/} in the literature.
Since, as we just argued, the expression~(\AA) is a special case of
the expression~(\AS), the Gau\ss ian polynomials are also subsumed
by expression~(\AS) as special cases. This explains Stanton's notion of
``fake Gaussian sequences."} While it seems that our Conjecture~\TE\
is a special case of the above conjecture, this is not so.
Stanton's conjecture requires polynomiality {\it for all\/}~$m$,
which is a big restriction.
We do not do that. However, on the basis of numerous computations
with ``random"~$m$ and ``random'' sequences\footnote{Aside from
ad hoc ``small number"-checks, we let several
computers run for several days checking sequences
$(0^s,a,b,b,a,0^s)$ and $(0^s,a,b,c,b,a,0^s)$,
with $a,b,c,s$ random integers between 0 and 10, and $m$
between $s$ and $1000+s$.} $(a_1,a_2,\dots,a_n)$ (see Remark~\RB.(2) after
the proof of Lemma~\TF\ for further evidence), we tend to believe
that Stanton's polynomiality restriction can be dropped and
our weaker hypothesis is sufficient.

\proclaim{Conjecture \TG}
Let $m$ and $n$ be given positive integers,
and let $(a_1,a_2,\dots,a_n)$ be a symmetric
sequence of non-negative integers.
If the expression~{\rm(\AS)}
is a polynomial\/ {\rm(}in~$q${\rm)} then it has non-negative coefficients.
\endproclaim

We furthermore believe that even symmetry is not essential here.

\proclaim{Conjecture \TH}
Let $m$ and $n$ be given positive integers with $m\ge n$,
and let $(a_1,a_2,\mathbreak\dots,a_n)$ be a
sequence of non-negative integers.
If the expression~{\rm(\AS)}
is a polynomial\/ {\rm(}in~$q${\rm)} then it has non-negative coefficients.
\endproclaim

The restriction of $m\ge n$ should be noted.
Stanton~\cite{\StatAD, bottom of p.~1} provides the example 
$m=1$ and
$(a_1,\dots,a_{17})=(1,3,1,1,1,1,1,1,2,1,1,1,1,1,1,1,1)$
where (\AS) is a polynomial that contains the term $-q^7$.
However, computer calculations up to $m=200$ indicate that
expression~(\AS) with the same $a_i$'s is a polynomial
(this was already shown by Stanton) with non-negative coefficients
not only for all $m\ge17$, but already for all $m\ge2$.

\medskip
(4)
We use our result for $l=2$ to confirm a conjecture of
Hopkins and Rubey~\cite{\HoRuAA} about a presumed cycling
sieving polynomial; see Corollary~\TC\ below.

\medskip
(5)
Conjecture~\TE\ is reminiscent of Warnaar and Zudilin's
``$q$-rious positivity"~\cite{\WaZuAB}. However, it is not
a special case since Landau's criterion is never satisfied
here (except in the trivial cases where $l=0$ or where $l=k$).
On the other hand, it is also not true that Warnaar and
Zudilin's conjecture can be generalised to a result where
one replaces Landau's criterion by the assumption of
polynomiality of the $q$-factorial expression. For example,
we have
$$\align
\frac {[12]_q!\,([2]_q!)^2} {[11]_q!\,[4]_q!\,[1]_q!}
&=1 - q^2 + q^3 + q^4 - q^5 + q^7,
\\
\frac {[10]_q!\,[4]_q!\,[3]_q!\,[1]_q!}
{[9]_q!\,[5]_q!\,([2]_q!)^2}
&=1 + q^2 - q^3 + q^4 + q^6,
\\
\frac {[12]_q!\,([2]_q!)^3}
{[11]_q!\,[4]_q!\,([1]_q!)^3}
&=1 + q - q^2 + 2 q^4 - q^6 + q^7 + q^8,
\\
\frac {[12]_q!\,([2]_q!)^5}
{[11]_q!\,[4]_q!\,([1]_q!)^7}
&=1 + 3 q + 2 q^2 - q^3 + q^4 + 4 q^5 + q^6 - q^7 + 2 q^8 + 3 q^9 + q^{10}.
\endalign
$$

\medskip
(6)
The expressions in (\AA) and (\AS) are (if polynomial) clearly
{\it cyclotomic generating functions} in the sense of~\cite{\BiSwAA}.
\endremark

\medskip
As announced, we now discuss the special cases where
$l=k,k-1,k-2,k-3$ and $l=0,1,2$.

\medskip
{\smc Case 1: $l=k$.}
This is trivial since the quotient in (\AA) reduces to~1.

\medskip
{\smc Case 2: $l=k-1$.}
In that case, the quotient in (\AA) reduces to 
$$
\frac {[n-k+1]_q} {[k]_q}
=\frac {1-q^{n-k+1}} {1-q^k}.
$$
Clearly, this is a polynomial if and only if $k$ divides $n-k+1$,
and then it is a terminating geometric series, which has non-negative
coefficients.

\medskip
{\smc Case 3: $l=k-2$.}
Here, the quotient in (\AA) reduces to 
$$
\frac {[n-k+2]_q\,[n-k+1]_q} {[k]_q\,[k-1]_q}
=\frac {(1-q^{n-k+2})(1-q^{n-k+1})} {(1-q^{k})(1-q^{k-1})}.
\tag\AL
$$
At this point, we should recall that
$$
q^n-1=\prod _{d\mid n} ^{}C_d(q),
\tag\AF
$$
where $C_d(q)$ denotes the $d$-th cyclotomic polynomial in $q$,
all of which are irreducible. Hence, the fraction~(\AL) can
only be a polynomial if both $k$ and $k-1$ divide one of
$n-k+1$ or $n-k+2$. If one of them divides $n-k+1$ and the
other divides $n-k+2$, then we are back in the situation
of Case~2, confirming the assertion of Conjecture~\TE.

It remains to discuss the case where both $k$ and $k-1$
divide the same number, that is, either $n-k+1$ or $n-k+2$.
In order to cope with this situation, we quote \cite{\KrMuAD, Cor.~6};
see also \cite{\ZudiAA, Ex.~1.8}.

\proclaim{Lemma \TD}
Let $a$ and $b$ be coprime positive integers, and let $\ga$ be an integer
with $\ga\ge (a-1)(b-1)$. Then the expression
$$
\frac {[ab]_q\,[\ga]_q} {[a]_q\,[b]_q}
$$
is a polynomial in $q$ with non-negative
integer coefficients.
\endproclaim

Now, if we suppose that both $k$ and $k-1$ divide $n-k+2$,
say $n-k+2=mk(k-1)$, then the fraction~(\AL) can be written in the
form 
$$
\frac {[k(k-1)]_q\,[n-k+1]_q} {[k]_q\,[k-1]_q}
\cdot\frac {[mk(k-1)]_q} {[k(k-1)]_q}.
$$
By Lemma~\TD, the first term in this product is a polynomial with
non-negative coefficients, while the second term is again a
terminating geometric series.

The case where both $k$ and $k-1$ divide $n-k+1$ is treated in the
same manner.

Altogether, this confirms Conjecture~\TE\ in the current case.

\medskip
{\smc Case 4: $l=k-3$.}
For this choice of $l$, the fraction (\AA) becomes
$$
\frac {[n-k+3]_q\,[n-k+2]_q\,[n-k+1]_q} {[k]_q\,[k-1]_q\,[k-2]_q}
=\frac {(1-q^{n-k+3})(1-q^{n-k+2})(1-q^{n-k+1})}
{(1-q^{k})(1-q^{k-1})(1-q^{k-2})}.
\tag\AM
$$

First of all, by the argument using cyclotomic polynomials that
we already used in Case~3,
for (\AM) to be a polynomial,
$k$, $k-1$, and $k-2$ must each divide one of
the arguments of the numerator factors, that is, one of
$n-k+3$, $n-k+2$, or $n-k+1$.
Following the arguments of Case~3, we must now distinguish between
several cases.

It might be that $k$, $k-1$, and $k-2$ each divide a different
factor\footnote{In abuse of language, we use abridged terminology here:
``factor" should actually be ``argument of a numerator factor".}
from the numerator. Then we are again back in the situation
of Case~2, confirming the assertion of Conjecture~\TE.

It might be that two of $k$, $k-1$, and $k-2$ divide the same
factor in the numerator, while the remaining one divides a
different factor. If, say, $k$ and $k-1$ divide the same factor
in the numerator, while $k-2$ divides a different one, then
Lemma~\TD\ applies again and shows that we obtain a polynomial
with non-negative coefficients. The other case is when
$k$ and $k-2$ divide the same factor in the numerator, while
$k-1$ divides a different one. This case has two subcases.
If $k$ and $k-2$ are relatively prime, then we argue as
before and see, using again Lemma~\TD, that we obtain a
polynomial with non-negative coefficients. However, it might be that
2 divides both $k$ and $k-2$. Counting cyclotomic factors on top
and in the bottom of~(\AM), we infer that there must also be
two terms out of $n-k+3$, $n-k+2$, and $n-k+1$ that are divisible
by~2. Clearly, these must be $n-k+3$ and $n-k+1$.
If $k-1$ divides $n-k+2$, then another application of Lemma~\TD\
leads to the desired conclusion of non-negativity of coefficients
of the polynomial~(\AM). If not, then the situation is more
complicated, see the next paragraph.

Continuing the previous case,
we suppose that $k$ and $k-2$ are even and divide $n-k+1$, and that
$k-1$ divides $n-k+3$. Let
$k=2K$ and $n-k+1=m\frac {k(k-2)} {2}=mK(2K-2)$ for appropriate positive
integers~$K$ and~$m$. Since $k-1$ divides $n-k+3$, we must have
$$
m\frac {k(k-2)} {2}+2
=2K(K-1)m+2\equiv0~\text{(mod $2K-1$)}.
$$
After simplification, we obtain the congruence
$$
(K-1)m\equiv-2~\text{(mod $2K-1$)}.
$$
Since 2 and $2K-1$ are relatively prime, we may multiply both sides of
the congruence by~2 and get
$$
2(K-1)m\equiv-4~\text{(mod $2K-1$)},
$$
or, after simplification,
$$
m\equiv4~\text{(mod $2K-1$)}.
$$
Consequently, we may write $m=4+M(2K-1)$, for some
non-negative integer~$M$. If we substitute these parametrisations
in~(\AM), then we obtain the fraction in~(\AN) below. In Lemma~\TF\
(also below)
it is shown that this fraction is a polynomial with non-negative coefficients.

In the other case, if $k$ and $k-2$ are even but divide $n-k+3$ while
$k-1$ divides $n-k+1$, one can proceed similarly. One finally arrives
at the fraction in~(\AO), and in Lemma~\TF\ it is shown that
this fraction is also a polynomial with non-negative coefficients.

Finally, we turn to the case where all of $k$, $k-1$, and $k-2$ divide
the same number out of $n-k+3$, $n-k+2$, and $n-k+1$. If $k$ is odd,
then $k$, $k-1$, and $k-2$ are pairwise coprime. Let us
write $A_1,A_2,A_3$ for $n-k+3,n-k+2,n-k+1$, in some order, such that
$k$, $k-1$, and $k-2$ divide~$A_1$. Using this notation,
we may rewrite the fraction in~(\AM) in the form
$$
\frac {[k(k-1)]_q\,[A_2]_q} {[k]_q\,[k-1]_q}\cdot
\frac {[k(k-1)(k-1)]_q\,[A_3]_q} {[k(k-1)]_q\,[k-2]_q}\cdot
\frac {[A_1]_q} {[k(k-1)(k-2)]_q}.
$$
Lemma~\TD\ applies to the first two terms in this product, while
the third term is a terminating geometric series.
This shows that the above fraction is a
polynomial with non-negative coefficients.

On the other hand, if $k$ is even then, arguing as before, also
$n-k+1$ and $n-k+3$ must be even. We use again the earlier notation,
with $k$, $k-1$, $k-2$ all dividing~$A_1$. In this case, we write
the fraction in~(\AM) in the form
$$
\frac {[k(k-2)/2]_q\,[A_2]_q} {[k]_q\,[k-2]_q}\cdot
\frac {[k(k-1)(k-2)/2]_q\,[A_3]_q} {[k(k-2)/2]_q\,[k-1]_q}\cdot
\frac {[A_1]_q} {[k(k-1)(k-2)/2]_q},
$$
where, here, $A_1$ and $A_2$ are chosen to be $n-k+1$ and $n-k+3$, in
the appropriate order. Lemma~\TD\ with $q$ replaced by~$q^2$ applies
to the first term, while Lemma~\TD\ applies directly to the second
term. The third term is again a terminating geometric series.

Altogether, this confirms Conjecture~\TE\ in our current case, Case~4, modulo the
subsequent lemma.

\proclaim{Lemma \TF}
Let $K$ and $M$ be non-negative integers.
Then the fractions
$$
\frac {
\prod _{i=0} ^{2}[4K(2K-2)+MK(2K-1)(2K-2)+i]_q}
{[2K]_q\,[2K-1]_q\,[2K-2]_q},\quad \text{for }K\ge2,
\tag\AN$$
and
$$
\frac {
\prod _{i=0} ^{2}[K(2K-5)(2K-2)+MK(2K-1)(2K-2)-i]_q}
{[2K]_q\,[2K-1]_q\,[2K-2]_q},\quad \text{for }K\ge3,
\tag\AO$$
are polynomials with non-negative coefficients.
\endproclaim

\demo{Proof}
We treat the first fraction. The second can be handled in a similar manner.

\medskip
First it should be observed that the factors in the numerator of~(\AN)
are explicitly
$$\multline
\left[2K(K-1)\big(4+M(2K-1)\big)\right]_q,\
\left[2K(K-1)\big(4+M(2K-1)\big)+1\right]_q,\\
\left[4K(2K-2)+MK(2K-1)(2K-2)+2\right]_q
=
\left[2(2K-1)\big((2K-1)+MK(K-1)\big)\right]_q,
\kern-1pt
\endmultline$$
so that, indeed, $2K$ and $2K-2$ divide the argument of the first factor,
while $2K-1$ divides the argument of the third factor.
Thus, all the cyclotomic polynomials $C_d(q)$ with $d\ge2$ dividing $1-q^{2K}$ 
divide the first factor, as well as
all the cyclotomic polynomials $C_d(q)$ with $d\ge2$ dividing $1-q^{2K-2}$.
There is one overlap though, namely $C_2(q)=1+q$ divides both $1-q^{2K}$
and $1-q^{2K-2}$, but only once the first factor. This is balanced by
the third factor, which is also divisible by $1+q$. Since
all the cyclotomic polynomials $C_d(q)$ with $d\ge2$ dividing $1-q^{2K-1}$
divide the third factor --- and $C_2(q)=1+q$ is not among them since
$2K-1$ is odd --- we have shown that the fraction~(\AN) is a
polynomial.

In order to make the following step more transparent, let us
temporarily put
$$\align
N_1&=(K-1)\big(4+M(2K-1)\big),\\
N_2&=4K+MK(2K-1),\\
N_3&=2\big((2K-1)+MK(K-1)\big).
\endalign$$
Using this notation, the fraction (\AN) can be written as
$$\multline
\frac {
[2KN_1]_q\,[(2K-2)N_2+1]_q\,[(2K-1)N_3]_q}
{[2K]_q\,[2K-1]_q\,[2K-2]_q}\\
=\left(1-q^{(2K-2)N_2+1}\right)\cdot
\frac {1-q^{2KN_1}} {1-q^{2K}}\cdot
\frac {1} {1-q^{2K-2}}\cdot
\frac {1-q^{(2K-1)N_3}} {1-q^{2K-1}}.
\endmultline$$
Expansion of geometric series turns this into
$$\multline
\sum_{j_1=0}^{N_1-1}\sum_{j_2\ge0}\sum_{j_3=0}^{N_3-1}
q^{2Kj_1+(2K-2)j_2+(2K-1)j_3}\\
-
\sum_{j_1=0}^{N_1-1}\sum_{j_2\ge0}\sum_{j_3=0}^{N_3-1}
q^{2Kj_1+(2K-2)j_2+(2K-1)j_3+(2K-2)N_2+1}.
\endmultline
\tag\AQ
$$
A subsum of the first sum is
$$\multline
\sum_{j_1=0}^{N_1-1}\sum_{j_4\ge0}\sum_{j_5=0}^{N_3-2}
q^{2Kj_1+(2K-2)(j_4+N_2-1)+(2K-1)(j_5+1)}\\
=
\sum_{j_1=0}^{N_1-1}\sum_{j_4\ge0}\sum_{j_5=0}^{N_3-2}
q^{2Kj_1+(2K-2)j_4+(2K-1)j_5+(2K-2)N_2+1}.
\endmultline$$
We see that this largely cancels with the second sum in~(\AQ), where only the terms
with $j_3=N_3-1$ remain. These form a series of order\footnote{The
order of a formal power series in~$q$ is the smallest exponent~$e$ such that
$q^e$ appears with non-zero coefficient in the series.}
$$
(2K-1)(N_3-1)+(2K-2)N_2+1
=
16K^2-18K+4+4MK(K-1)(2K-1).
\tag\AP$$
The total degree of the polynomial (as we now know) in~(\AN) is
$$
24K(K-1)+6MK(K-1)(2K-1)-6K+6.
$$
The order in (\AP) is larger than half of this degree.
Hence, by~(\AQ) and the above considerations,
we have shown that the polynomial in~(\AN) has non-negative coefficients
at least up to terms of half of the degree of the polynomial. Since
the polynomial in~(\AN) is a reciprocal polynomial\footnote{A polynomial
$P(q)$ in $q$ of degree $m$ is called {\it reciprocal\/} if
$P(q)=q^mP(1/q)$.} it follows that also the other coefficients must be
non-negative.\quad \quad \qed
\enddemo

\remark{Remark \RB}
(1)
Clearly, we could continue in this manner and next consider the case
where $l=k-4$, etc. However, it is equally clear that the
corresponding arguments get more and more involved. Since we did not
see how to build a general reasoning that would work for any $l$
and~$k$, we stop these considerations at this point and turn
to cases where the value of~$l$ is small.

\medskip
(2) Equipped with the previous arguments, we may also verify
further special cases of Conjecture~\TG.

If the sequence $(a_1,a_2,\dots,a_n)$ is of the form $(0^s,a,0^s)$ for
some non-negative integers $a$ and~$s$, then, by taking the $a$-th root
of expression~(\AS), we see that we are back in Case~2. Similarly,
if the sequence $(a_1,a_2,\dots,a_n)$ is of the form $(0^s,a,a,0^s)$,
then we are back in Case~3. 

If the sequence $(a_1,a_2,\dots,a_n)$ is of the form $(0^s,a,b,a,0^s)$,
then we see that the arguments in Case~4 can be almost straightforwardly
adapted to cover this slightly more general case as well.
\endremark

\medskip
{\smc Case 5: $l=0$.} Then the quotient in (\AA) reduces to the
$q$-binomial coefficient $\left[\smallmatrix n\\k\endsmallmatrix\right]_q$.
It is well-known that this has non-negative coefficients, see e.g\. \cite{\AndrAF,
Theorem~3.2}.

\medskip
{\smc Case 6: $l=1$.} Then the quotient in (\AA) reduces to 
$\frac {1} {[n]_q}\left[\smallmatrix n\\k\endsmallmatrix\right]_q$.
If $\gcd(k,n-k)=1$, this is known as a {\it rational
$q$-Catalan polynomial}, see e.g.~\cite{\ArmDAB, Sec.~5}. 
It is known that it has non-negative
coefficients. However, since we were not able to find a proof
of the {\it equivalence} that 
$\frac {1} {[n]_q}\left[\smallmatrix n\\k\endsmallmatrix\right]_q$
is a polynomial in~$q$
{\it if and only if} $\gcd(k,n-k)=1$ in the literature, for completeness
and for the
convenience of the reader we provide a proof of all these facts.

\proclaim{Theorem \TA}
The rational function 
$$\frac {1} {[n]_q}\left[\matrix n\\k\endmatrix\right]_q
\tag\AE$$
is a polynomial if and only if $\gcd(k,n-k)=1$. In the latter
case, it is a polynomial with non-negative integer coefficients.
\endproclaim

\demo{Proof}
By (\AF), the rational function in~(\AE) can be decomposed into
a product of cyclotomic polynomials (with positive and negative
exponents). The cyclotomic
polynomial $C_d(q)$ appears in~(\AE) with exponent
$$
\fl{\tfrac {n-1} {d}}-\fl{\tfrac {k} {d}}
-\fl{\tfrac {n-k} {d}}+\chi(d=1).
\tag\AH$$
Here, $\chi(\Cal A)=1$ if $\Cal A$ is true and $\chi(\Cal A)=0$ otherwise.
We must show that this is non-negative for all~$d\ge1$ if and only if
$\gcd(k,n-k)=1$.

We first note that (\AH) vanishes for  $d=1$. From now on, let $d\ge2$.
Writing $N=\left\{\frac {n-1} {d}\right\}$ and
$K=\left\{\frac {k} {d}\right\}$, where $\{\al\}$ denotes the
fractional part of a real number~$\al$, we may convert (\AH) into
$$
\fl{N}-\fl{K}
-\fl{N-K+\tfrac {1} {d}}
=-\fl{N-K+\tfrac {1} {d}}.
\tag\AG$$
This is non-negative except when $N-K+\tfrac {1} {d}=1$.
That equation however is equivalent with $N=\frac {d-1} {d}$
and $K=0$. This means that $d\mid n$ and $d\mid k$, which implies
$\gcd(k,n-k)>1$.
This proves the first assertion in the statement of the theorem.

\medskip
Let $\gcd(k,n-k)=1$.
For the proof of the second assertion, we reproduce an argument that
can be found in \cite{\AndrCB, Theorem~2} or \cite{\ReSWAB, Prop.~10.1(iii)}.
We recall that it is
well-known that the $q$-binomial coefficient
$\left[\smallmatrix n\\k\endsmallmatrix\right]_q$ is a
reciprocal and unimodal\footnote{A polynomial $P(q)=
\sum _{i=0} ^{m}p_iq^i$ is called {\it unimodal\/} if there is an
integer $r$ with $0\le r\le m$ and $0\le p_0\le\dots\le p_r\ge\dots\ge
p_m\ge0$. It is well-known that $q$-binomial coefficients are
unimodal; see \cite{\StanBI, Ex.~7.75.d}.} polynomial in~$q$
of degree~$k(n-k)$. As a consequence, the expression
$$
(1-q)\left[\matrix n\\k\endmatrix\right]_q
$$
is a polynomial with a non-negative coefficient of $q^s$ for
$0\le s\le \frac {1} {2}k(n-k)$. We multiply this expression
by the geometric series $\frac {1} {1-q^n}=1+q^n+q^{2n}+\cdots$.
It follows that
$$
\frac {1-q} {1-q^n}\left[\matrix n\\k\endmatrix\right]_q
=\frac {1} {[n]_q}\left[\matrix n\\k\endmatrix\right]_q
\tag\AI$$
is a polynomial (by the first part, since we assumed that
$\gcd(k,n-k)$=1) with a non-negative coefficient of~$q^s$
for $0\le s\le \frac {1} {2}k(n-k)$. Since also (\AI) is a
reciprocal polynomial, all the other coefficients are non-negative
as well.\quad \quad \qed
\enddemo

\medskip
{\smc Case 7: $l=2$.} Then the quotient in (\AA) reduces to 
$$\frac {[1]_q\,[2]_q} {[n]_q\,[n-1]_q}\bmatrix n\\k\endbmatrix_q
=\frac {(1-q)(1-q^2)} {(1-q^n)(1-q^{n-1})}\bmatrix n\\k\endbmatrix_q.
\tag\AB
$$
We show in the theorem below that, if this is a polynomial,
it has non-negative coefficients.

\proclaim{Theorem \TB}
The rational function {\rm(\AB)}
is a polynomial if and only if
$$
\{\gcd(k,n-k)
,\gcd(k,n-k-1)
,\gcd(k-1,n-k)\}\subseteq\{1,2\}.
\tag\AD
$$
In that
case, it is a polynomial with non-negative integer coefficients.
\endproclaim

\demo{Proof}
Let again $C_d(q)$ denote the $d$-th cyclotomic polynomial in $q$.
By~(\AF), the cyclotomic
polynomial $C_d(q)$ appears in~(\AB) with exponent
$$
\fl{\tfrac {n-2} {d}}-\fl{\tfrac {k} {d}}
-\fl{\tfrac {n-k} {d}}
+2\chi(d=1)+\chi(d=2).
\tag\AJ$$
We must show that this is non-negative for all~$d\ge1$ if and only if
the condition in~(\AD) holds.

We begin with $d=1$: in that case (\AJ) equals zero.

Next let $d=2$. In that case, the expression (\AJ) reduces to
$$
\fl{\tfrac {n} {2}}-\fl{\tfrac {k} {2}}
-\fl{\tfrac {n-k} {2}}.
$$
It is easy to see that this is always non-negative.

From now on, let $d\ge3$.
Writing $N=\left\{\frac {n-2} {d}\right\}$ and
$K=\left\{\frac {k} {d}\right\}$, we may convert (\AJ) into
$$
\fl{N}-\fl{K}
-\fl{N-K+\tfrac {2} {d}}
=-\fl{N-K+\tfrac {2} {d}}.
\tag\AK$$
This is non-negative except if $N-K+\frac {2} {d}=1$
or $N-K+\frac {2} {d}=1+\frac {1} {d}$. These exceptions translate
into the three cases
$$
(N,K)\in\left\{(\tfrac {d-2} {d},0),\,
(\tfrac {d-1} {d},\tfrac {1} {d}),\,
(\tfrac {d-1} {d},0)\right\}.
$$
If $(N,K)=(\tfrac {d-2} {d},0)$ then $d\mid n$ and $d\mid k$,
implying $d\mid\gcd(k,n-k)$. This violates~(\AD).
If $(N,K)=(\tfrac {d-1} {d},\frac {1} {d})$ then $d\mid (n-1)$
and $d\mid (k-1)$,
implying $d\mid\gcd(k-1,n-k)$. Again, this violates~(\AD).
Finally, if $(N,K)=(\tfrac {d-1} {d},0)$ then $d\mid (n-1)$ and $d\mid k$,
implying $d\mid\gcd(k,n-k-1)$. This is also a violation of~(\AD).

We have proved the first assertion in the statement of the theorem.

\medskip
Now let the condition in (\AD) be satisfied.
For the proof of the second assertion, we start with
$$
\frac {1-q} {1-q^n}\bmatrix n\\k\endbmatrix_q\quad \text{respectively}
\quad 
\frac {1-q} {1-q^{n-1}}\bmatrix n-1\\k-1\endbmatrix_q,
\tag\AC
$$
depending on whether $\gcd(k,n-k)=1$ or $\gcd(k,n-k)=2$.
Given~(\AD), 
both expressions in~(\AC) are rational
$q$-Catalan polynomials, which have even degree, as is easily seen.
Furthermore, it follows from \cite{\StucAA, Theorem~1.2 with
$k=0$} that they are parity-unimodal\footnote{A polynomial $P(q)=
\sum _{i=0} ^{m}p_iq^i$ is called {\it parity-unimodal\/} if the
``even polynomial" $\sum _{i=0} ^{\fl{m/2}}p_{2i}q^i$
and the ``odd polynomial" $\sum _{i=0} ^{\fl{(m-1)/2}}p_{2i+1}q^i$
are both unimodal.}. Combined with the fact that they are also
reciprocal polynomials, these observations imply that
the products
$$
(1-q^2)\times\frac {1-q} {1-q^n}\bmatrix n\\k\endbmatrix_q\quad
\text{respectively}\quad
(1-q^2)\times\frac {1-q} {1-q^{n-1}}\bmatrix n-1\\k-1\endbmatrix_q
$$
have non-negative coefficients of $q^s$ for
$s\le \frac {1} {2}\big(k(n-k)-(n-1)\big)$
respectively for $s\le \frac {1} {2}\big((k-1)(n-k)-(n-2)\big)$.
If we multiply the first expression by the
geometric series $\frac {1} {1-q^{n-1}}$, respectively the second
expression by $\frac {1} {1-q^k}$,
we obtain the series
$$
\frac {1-q^2} {1-q^{n-1}}
\times\frac {1-q} {1-q^n}\bmatrix n\\k\endbmatrix_q\quad
\text{respectively}\quad 
\frac {1-q^2} {1-q^k}
\times\frac {1-q} {1-q^{n-1}}\bmatrix n-1\\k-1\endbmatrix_q
$$
which have also non-negative coefficients of $q^s$ for
$s\le \frac {1} {2}\big(k(n-k)-(n-1)\big)$ respectively
for $s\le \frac {1} {2}\big((k-1)(n-k)-(n-2)\big)$.
However, both expressions equal the expression~(\AB), of which we know
that it is a polynomial. It is moreover a reciprocal polynomial,
therefore also all other coefficients are non-negative.\quad \quad \qed
\enddemo

As an application, we show that
$$
\frac {
\prod _{j=1} ^{3n}(1-q^{2j})}
{
\prod _{j=2} ^{2n+1}(1-q^j)
\prod _{j=2} ^{n+1}(1-q^{2j})}
\tag\AR$$
is a polynomial with non-negative integer coefficients.
This was conjectured by Hopkins and Rubey
in \cite{\HoRuAA, Conj.~6.2} as part of
a presumed cyclic sieving phenomenon for Kreweras words.

\proclaim{Corollary \TC}
The expression {\rm(\AR)} is a polynomial with non-negative integer
coefficients.
\endproclaim

\demo{Proof}
For convenience, we introduce the usual notation for $q$-shifted
factorials\break $(\al;q)_m:=(1-\al)(1-\al q)\cdots(1-\al q^{m-1})$ for
$m\ge1$, and $(\al;q)_0:=1$.
Then (\AR) can be rewritten as
$$
\frac {(q^2;q^2)_{3n}} {(q^2;q)_{2n}\,(q^4;q^2)_n}
=
\frac {(q^2;q^2)_{3n}\,(-q^2;q)_{2n}} {(q^4;q^2)_{2n}\,(q^4;q^2)_n}
=
\frac {(1-q^2)(1-q^4)\,(-q^3;q)_{2n-1}}
{(1-q^{6n+4})(1-q^{6n+2})}
\bmatrix 3n+2\\n+1\endbmatrix_{q^2}.
$$
Using Theorem~\TB\ with $q$ replaced  by~$q^2$,
$n$ replaced by $3n+2$ and $k=n+1$, and the fact that
$(-q^3;q)_{2n-1}$ is a polynomial in~$q$ with non-negative
integer coefficients, we obtain the assertion of the
corollary.\quad \quad \qed
\enddemo

\remark{Remark}
The starting point for this work was indeed
Conjecture~6.2 by Hopkins and Rubey~\cite{\HoRuAA}.
Our investigations around this conjecture led us first to
extract what has now become Theorem~\TB\ as the essential
germ. Later on we discovered the more general conjecture
on quotients of
$q$-binomial coefficients that is in the centre of this note.
\endremark

\remark{Acknowledgement}
We thank Wadim Zudilin for inspiring discussions, and in particular
for suggesting to also consider the cases of Conjecture~\TE\
where the difference between $k$ and~$l$ is fixed.
Furthermore, we are indebted to Sam Hopkins for pointing out~\cite{\StatAD}.
\endremark


\Refs

\ref\no \AndrAF\by G. E. Andrews \yr 1976 \book The Theory of
Partitions\publ Addison--Wesley\publaddr Reading \bookinfo
Encyclopedia of Math\. and its Applications, Vol.~2\finalinfo
reprinted by Cambridge University Press, Cambridge, 1998\endref

\ref\no \AndrCB\by G. E. Andrews \yr 2004 \paper The 
Friedman--Joichi--Stanton monotonicity conjecture at primes\inbook 
Unusual Applications of Number Theory\ed M.~Nathanson\publ 
DIMACS Ser\. Discrete Math\. Theor\. Comp\. Sci., vol.~64, 
Amer\. Math\. Soc\.\publaddr Providence, R.I.\pages 9--15\endref

\ref\no \ArmDAB\by D. Armstrong\paper Lattice points and rational
$q$-Catalan numbers\jour preprint, {\tt ar$\chi$iv:2403.06318}
\vol  \pages \endref

\ref\no \BiSwAA\by S. C. Billey and J. P. Swanson\paper Cyclotomic
generating functions\jour Electron\. J. Combin\. \vol 31\rm(4)
\yr 2024\pages Paper No.~4.4, 43~pp\endref



\ref\no\HoRuAA \by S. Hopkins and M. Rubey\paper Promotion of Kreweras words
\jour Selecta Math. (N.S.)\vol 28\yr 2022\pages Paper No.~10, 38~pp\endref

\ref\no \KrMuAD\by C. Krattenthaler and T. W. M\"uller \yr 2013\paper
Cyclic sieving for generalised 
non-crossing partitions associated with complex reflection groups of
exceptional type 
\inbook Advances in Combinatorics, volume in memory of Herbert S.~Wilf,
\eds I.~Kotsireas, E.~Zima\publ Springer-Verlag \pages 209--247\endref

\ref\no \ReSWAB\by V.    Reiner, D. Stanton and D. White \yr 2004
\paper The cyclic sieving phenomenon\jour J. Combin\. Theory
Ser.~A\vol 108\pages 17--50\endref 

\ref\no \StanBI\by R. P. Stanley \yr 1999 \book Enumerative Combinatorics
\bookinfo Vol.~2\publ Cambridge University Press\publaddr Cambridge\endref

\ref\no \StatAD\by D.    Stanton \yr \paper Fake Gaussian sequences
\jour unpublished manuscript, 6~pages, 1998; available at\linebreak
{\tt https://www-users.cse.umn.edu/\~{}stant001/PAPERS/macmahon.pdf}\vol \pages \endref

\ref\no\StucAA \by E. N. Stucky\paper
Parity-unimodality and a cyclic sieving phenomenon for necklaces
\jour SIAM J. Discrete Math\.\vol 35\yr2021\pages 2049--2069\endref

\ref\no\WaZuAB \by S. O. Warnaar and W. Zudilin\paper
A $q$-rious positivity\jour
Aequationes Math\.\vol 81\yr 2011\pages 177--183\endref

\ref\no \ZudiAA\by W. Zudilin\book
Analytic Methods in Number Theory --— When Complex Numbers Count
\publ Monogr. Number Theory, vol.~11\publaddr
World Scientific Publishing Co\. Pte\. Ltd., Hackensack, NJ\yr 2024\endref

\endRefs
\enddocument